\documentclass[12pt]{article}
\usepackage{mathtools} \mathtoolsset{showonlyrefs} 
\usepackage{amsfonts}
\usepackage[utf8]{inputenc}
\usepackage{amsmath,amssymb,amsthm}
\usepackage{mathdots}
\usepackage{graphicx}

\usepackage[colorlinks=false,urlcolor=blue,pdfpagemode=None]{hyperref}

\usepackage{mathrsfs}
\usepackage{bm}
\usepackage[margin=0.9in]{geometry}

\newtheorem{thm}{Theorem}[section]
\newtheorem{prop}[thm]{Proposition}
\newtheorem{lemma}[thm]{Lemma}
\newtheorem{cor}[thm]{Corollary}
\theoremstyle{definition}
\newtheorem{definition}[thm]{Definition}
\newtheorem{remark}[thm]{Remark}

\usepackage{tikz} 

\usepackage{mdframed}
\mdfdefinestyle{MyFrame}{%
    linecolor=black,
    outerlinewidth=3pt,
    innertopmargin=\baselineskip,
    innerbottommargin=\baselineskip,
    innerrightmargin=10pt,
    innerleftmargin=10pt,
    backgroundcolor=gray!10!white}

\newcommand{\calA}{\mathcal{A}}
\newcommand{\calB}{\mathcal{B}}
\newcommand{\calC}{\mathcal{C}}
\newcommand{\calD}{\mathcal{D}}

\newcommand{\calG}{\mathcal{G}}
\newcommand{\calH}{\mathcal{H}}

\newcommand{\bbD}{\mathbb{D}}
\newcommand{\bbE}{\mathbb{E}}

\newcommand{\bbG}{\mathbb{G}}

\newcommand{\bbL}{\mathbb{L}}

\newcommand{\bbP}{\mathbb{P}}

\newcommand{\bbR}{\mathbb{R}}

\newcommand{\bbV}{\mathbb{V}}

\newcommand{\bbZ}{\mathbb{Z}}

\newcommand{\sfC}{{\sf C}}
\newcommand{\sfN}{{\sf N}}

\newcommand{\bfS}{{\bf S}}
\newcommand{\bfT}{{\bf T}}
\newcommand{\bfZ}{{\bf Z}}
\newcommand{\La}{\Lambda}
\newcommand{\eps}{\varepsilon}

\newcommand{\comment}[1]{}

\title{Universality of Large-Scale Geometry For Critical Random-Cluster Models}
\author{Ulrik Thinggaard Hansen\footnote{Institute for Stochastics, KIT, ulrik.hansen@kit.edu }  \text{ } and Ioan Manolescu\footnote{Department of Mathematics, UniFr, ioan.manolescu@unifr.ch}}

\begin{document}
\maketitle
\begin{abstract}
The understanding of the behaviour of critical lattice models in two dimensions has made immense progress in the past couple of decades, with conformal invariance as its conjectured end goal. 
We consider here the critical random-cluster model (or FK-percolation) on a large class of graphs --- namely bi-periodic isoradial graphs --- with cluster-weight~$q \in [1,4]$.
Following up on the establishment of rotational invariance at criticality in \cite{Rotation}, we prove that the large scale behaviour of these models is asymptotically the same on all such graphs.
These include the critical homogeneous models on the square, triangular and hexagonal lattices. 
\end{abstract}

\section{Introduction}

\subsection{Motivation}
The successful import of the renormalisation group formalism from quantum field theory was a landmark achievement in statistical mechanics\footnote{We refer to \cite{RGHistory} for a historical exposition.} that, while remaining elusive from a rigorous point of view, yielded two concrete predictions:
\begin{enumerate}
\item[1.] Critical statistical mechanics models should be universal in the sense that their scaling limits should not depend on the choice of underlying lattice.
\item[2.] The scaling limits should enjoy conformal symmetries.
\end{enumerate}
 In particular, since the planar conformal group is infinite dimensional, scaling limits of 2D statistical mechanics models should enjoy incredibly rich symmetries.
 
 With Schramm's introduction of the SLE in \cite{SLEBirth}, this conjectural picture took on a much more tangible form, and the subsequent twenty years have seen a staggering amount of progress, establishing the conjecture of conformal invariance in a number of different models (see \cite{Ave+26,ADCHLM26, CamiaNewmanFull, IsingInvariance, Giuliani1, Giuliani2, Dominos, HarmExp, GFF, TriLat} for a non-exhaustive list).
 
 However, with the exception of \cite{ADCHLM26,Giuliani1,Giuliani2}, these approaches are characterised by exploiting very particular features of discrete holomorphicity in the models and as such, do not offer robust methods for treating other models. 
 
Another approach is to focus on universality. That is, rather than taking a single model and trying to compute its scaling limit directly, take a class of related models and show that all of their possible scaling limits must agree. 
 
 Such an approach has also seen some development in the last fifteen years, starting with bond percolation in \cite{ManoGrimmett}, extending to more general random-cluster models in \cite{JiHuangLi} and lately, being used to establish not only universality across a class of square lattices, but rotational invariance of any possible scaling limit in \cite{Rotation}.
 
 In this paper, we combine the results of \cite{Rotation} with some of the conceptual techniques from \cite{ManoGrimmett} in order to extend the universality result of \cite{Rotation}  to a larger class of graphs. As a special case, we establish large scale rotational invariance of the critical random-cluster model on any bi-periodic isoradial lattice.

\subsection{Main Result}

In order to state our results properly, we give a brief introduction to some of the objects with which we will be working, all of which will be treated in detail in Section~\ref{sec:background}.

We say that an embedded planar graph~$\mathbb{G}=(\mathbb{V},\mathbb{E})$ is \textbf{isoradial} if every face of~$\mathbb{G}$ is inscribed in a circle of radius~$1$ and the centre of the corresponding circle lies within that face. 
We say that~$\mathbb{G}$ is \textbf{bi-periodic} if there exist two linearly independent vectors~$v,w\in \mathbb{Z}^2$ such that~$\mathbb{G}$ is invariant under translation by~$v$ and~$w$. 

For a finite graph~$G=(V,E)$, the \textbf{random-cluster model} is the probability measure on~$\{0,1\}^E$ given by
~$$
\phi_{\boldsymbol{p},q,G}[\{\omega\}]=\frac{1}{Z_{\boldsymbol{p},q,G}}q^{\kappa(\omega)}\prod_{e\in E} p_e^{\omega(e)} (1-p_e)^{1-\omega(e)},
$$
 where~$\boldsymbol{p}=(p_e)_{e\in E}\in (0,1)^E$ is a vector of probabilities,~$q\in (0,\infty)$ is called the cluster-weight,~$\kappa$ is the number of components in~$\omega$ and~$Z$ is a normalising constant called the \textbf{partition function}.

To each edge~$e$ of an isoradial graph, we shall attach a certain weight~$p_e$ depending only on the local geometry of the graph in a way that shall be made precise later on. For this choice of weights, we shall simply refer to the corresponding finite volume measures as~$\phi_{G}$, the so-called critical model with isoradial weights. By usual weak compactness elements, this gives rise to corresponding infinite-volume measures~$\phi_{\mathbb G}$ on infinite isoradial graphs~$\mathbb{G}$, which turn out to be unique in the regime of interest for us, namely~$q \in [1,4]$.
We will also write ~$\phi_{\delta \mathbb G}$ for this same measure on the graph~$\mathbb G$ rescaled by~$\delta>0$.

An example of an isoradial graph is the square lattice~$\sqrt 2 \bbZ^2$; note  the rescaling factor~$\sqrt 2$ that makes the faces of~$\bbZ^2$ be inscribed in circles of radius~$1$ rather than~$1/\sqrt 2$. 
Then,~$\phi_{\sqrt 2 \bbZ^2}$ is the critical random-cluster measure \cite{BefDum12} with edge-parameters~$\frac{\sqrt q}{1 + \sqrt q}$ for all edges, which has been extensively studied in recent years \cite{FKDiscont,ScalingRelations, SuperRSW, DCCont,GlazmanLamers,RaySpinka}. 

It is among this class of random-cluster models that universality is expected when~$q\in [1,4]$ -- that is, when the model on~$\mathbb{Z}^2$ undergoes a continuous phase transition (see \cite{FKDiscont,DCCont} with more recent proofs in \cite{RaySpinka,GlazmanLamers}). That the phase transition is of the same type on all bi-periodic isoradial graphs was proved in \cite{JiHuangLi}.
 
To make a precise statement, we fix a metric that is naturally associated with re-scaled percolation configurations. 
For two oriented loops,~$\gamma,\gamma':S^1\to \mathbb{C}$, we define their natural distance
\begin{align}\label{eq:single_loop_dist}
d(\gamma,\gamma')=\inf_{\varphi}\| \gamma\circ \varphi-\gamma'\|_{\infty},
\end{align}
where~$\varphi$ ranges over all orientation-preserving homeomorphisms of~$S^1$ into itself.

For two sets of oriented loops~$\mathcal{F}$ and~$\mathcal{G}$, 
we define their \textbf{Camia-Newman distance}~$d_{\rm CN}$ such that
$
d_{\rm CN}(\mathcal{F},\mathcal{G})\leq \varepsilon$ if and only if for every loop~$\gamma\in \mathcal{F}$ with~$\|\gamma\|_{\infty}\leq \frac{1}{\varepsilon}$, 
there exists a loop~$\gamma'\in \mathcal{G}$ such that~$d(\gamma,\gamma')<\varepsilon$, and vice-versa\footnote{It is worth noting that the Camia-Newman distance is essentially only a function of loops which are big. Indeed, $d_{\rm CN}(\gamma,\gamma')$ is bounded above by the sums of the diameters of the loops and the Hausdorff distance between the traces. As such, small loops will not contribute as long as there are any loops of similar diameter nearby in the other configuration.}. 

For convenience, we will also sometimes consider the global\footnote{Note that the global metric takes values in $[0,\infty]$ but otherwise satisfies the usual axioms of a metric. We will only use it on pairs of configurations where the distance is finite.} Camia-Newman metric~${d}_{\rm global\, CN}$ which is defined similarly, except that we do not restrict our attention to loops with~$\|\gamma\|_{\infty}\leq \frac{1}{\varepsilon}$.

Finally, for two probability measures~$\phi, \phi'$ on oriented families of loops, write 
\begin{align}\label{eq:dis_law}
d_{\rm CN}(\phi,\phi') := \inf\big\{\varepsilon>0 \mid \exists \text{ coupling~$\bbP$ of~$\phi$ and~$\phi'$ with }	\mathbb P[d_{\rm CN}(\omega,\omega')> \varepsilon] <\varepsilon\big\}.
\end{align}

In our applications, the loop configurations will be the set of cluster interfaces surrounding primal and dual clusters of a percolation  configuration~$\omega$ on some isoradial graph~$\bbG$. The exact definition is given in Section~\ref{sec:background}, see also Figure~\ref{fig:loop_rep}. For now, the only thing to mention is that these loops will always be oriented so that the dual edges lie on their left. 

The identification between percolation and loop configurations is one-to-one, and we shall use the Camia-Newman distance between percolation configurations by extension.

\begin{thm}\label{thm:main}
	Fix~$q \in [1,4]$.
	For any bi-periodic isoradial graph~$\mathbb{G}$, there exists~$c >0$ such that, for all~$\delta \leq 1$, 
	\begin{align}\label{eq:main_result}
		d_{\rm CN}\big(\phi_{\delta \bbG} , \phi_{\delta \sqrt 2 \bbZ^2} \big) \leq \delta^c.
	\end{align}
\end{thm}

The assumption of bi-periodicity of $\bbG$ above may be weakened, in the spirit of \cite{ManoGrimmett}, but not completely removed. 
We believe the result to hold for general isoradial graphs, but our proof does rely on bi-periodicity. 

Two issues arise when working with non-periodic graphs: the potential existence of edges of length arbitrarily close to $0$ or $2$
(this is excluded by the {\em bounded angles property} of  \cite{ManoGrimmett})
and the potential lack of regularity within $\bbG$ (such as the existence of a square grid structure, which was termed the {\em square grid property} in  \cite{ManoGrimmett}).
To extend Theorem~\ref{thm:main} to general isoradial graphs, one would need to tackle these difficulties. 
The authors believe they are best tackled together and that they are surmountable, but do require new ideas. 

\begin{cor}
For any~$q \in [1,4]$ and any bi-periodic isoradial graph~$\mathbb{G},$ any sub-sequential limit of~$\phi_{\delta \mathbb{G}}$ as~$\delta \searrow 0$  is invariant under rotations. 
\end{cor}

\begin{proof}
Consider a bi-periodic, isoradial graph~$\bbG$ and a rotation~$\widetilde \bbG$ of this same graph. 
As they are both bi-periodic and isoradial, \eqref{eq:main_result} implies 
\begin{align}
	d_{\rm CN}\big(\phi_{\delta \bbG} ,\phi_{\delta \widetilde \bbG} \big) \xrightarrow[\delta \to 0]{}0.
\end{align}
This implies the result. 
\end{proof}

It was recently proved in \cite{Ave+26} that the critical random-cluster model on $\bbZ^2$ with cluster weight $q = 4$ converges to the conformal loop ensemble CLE$_4$. 
By Theorem~\ref{thm:main}, this extends readily to $\phi_{\bbG}$ for any bi-periodic isoradial graph $\bbG$.

\begin{cor}\label{cor:CLE}
Fix~$q = 4$ and~$\mathbb{G}$ a  bi-periodic isoradial graph. Then~$\phi_{\delta \mathbb{G}}$ converges, as~$\delta \searrow 0$, to the full-plane nested CLE$_4$. 
This includes the homogeneous models on the triangular and hexagonal lattices with critical parameters $1/2$ and $4/5$, respectively. 
\end{cor}

The full-plane CLE$_4$ is the local limit of the process in increasing domains; see \cite{KemWer16}. The convergence above holds for the $d_{\rm CN}$ distance. 

We do not claim the convergence in domains of $\bbG$. Indeed, to deduce the convergence in domains from that in the full plane, one requires a form of continuity of the measure in the shape of the domain in which it is defined. Such a continuity for $\phi_{\bbZ^2}$ is the byproduct of the convergence in \cite{Ave+26}. For isoradial graphs, we do not have an a priori proof of such  continuity, which is why we restrict Corollary~\ref{cor:CLE} to the full plane. 

For $q = 2$, the scaling limit of the random-cluster model (also known in this case as the FK-Ising model) was determined on isoradial graphs and beyond in \cite{CheSmi12, Che20} using direct analysis of each model. For $q \in [1,2)\cup (2,4)$, the scaling limit on $\bbZ^2$ and isoradial graphs alike remains conjectural.

\subsubsection*{Structure of the Paper}
Section~\ref{sec:background} is dedicated to the necessary background material for the paper. In particular, we review isoradial graphs and the associated critical random-cluster models as well the star-triangle transform and its consequences. The latter is our main tool for comparing a priori different models.

Section~\ref{sec:cube_flip} then focuses on a particular composite type of transformation, which we will call a cube flip; it is essentially a large-scale star-triangle transformation. It is reminiscent of the interchanging of square lattices in \cite{Rotation}, but in finite volume.  
The estimates of this section are crucial, but follow  arguments similar to~\cite{Rotation}.
To avoid repeating tedious details, proofs here will be sketched, but may be followed without prior knowledge of \cite{Rotation}.

Section~\ref{sec:main_proof} turns to the proof of Theorem~\ref{thm:main}. First, the graph is transformed so as to be composed of large cubes; these transformations are mesoscopic, and therefore incur an error which is easily controlled. Then, rapid fire applications of the cube flip of Section~\ref{sec:cube_flip} transform this graph into a square lattice, and the estimates obtained there easily bound the ensuing errors. 

\section{Background}\label{sec:background}

\subsection{Isoradial Graphs}\label{sec:isoradial_graphs}

\begin{figure}
\begin{center}
\includegraphics[scale=0.8]{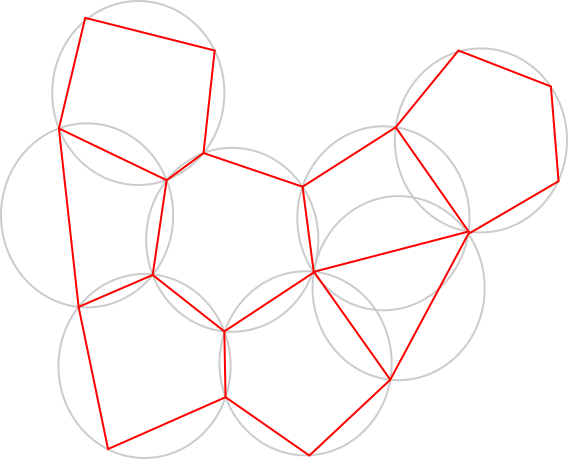} 
\includegraphics[scale=0.8]{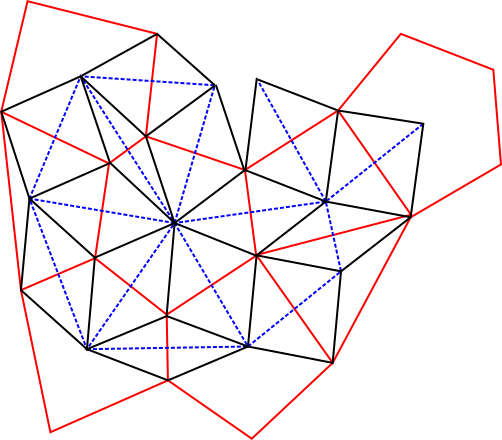} \scalebox{.8}{
\begingroup%
  \makeatletter%
  \providecommand\color[2][]{%
    \errmessage{(Inkscape) Color is used for the text in Inkscape, but the package 'color.sty' is not loaded}%
    \renewcommand\color[2][]{}%
  }%
  \providecommand\transparent[1]{%
    \errmessage{(Inkscape) Transparency is used (non-zero) for the text in Inkscape, but the package 'transparent.sty' is not loaded}%
    \renewcommand\transparent[1]{}%
  }%
  \providecommand\rotatebox[2]{#2}%
  \newcommand*\fsize{\dimexpr\f@size pt\relax}%
  \newcommand*\lineheight[1]{\fontsize{\fsize}{#1\fsize}\selectfont}%
  \ifx\svgwidth\undefined%
    \setlength{\unitlength}{164.40944882bp}%
    \ifx\svgscale\undefined%
      \relax%
    \else%
      \setlength{\unitlength}{\unitlength * \real{\svgscale}}%
    \fi%
  \else%
    \setlength{\unitlength}{\svgwidth}%
  \fi%
  \global\let\svgwidth\undefined%
  \global\let\svgscale\undefined%
  \makeatother%
  \begin{picture}(1,0.82758621)%
    \lineheight{1}%
    \setlength\tabcolsep{0pt}%
    \put(0,0){\includegraphics[width=\unitlength,page=1]{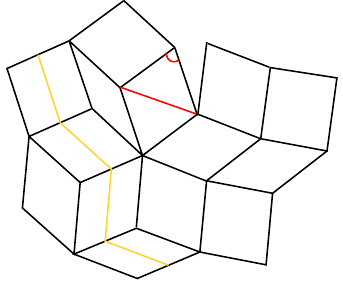}}%
    \put(0.42155861,0.48675557){\makebox(0,0)[lt]{\lineheight{1.25}\smash{\begin{tabular}[t]{l} $e$\end{tabular}}}}%
    \put(0.4731714,0.58539866){\makebox(0,0)[lt]{\lineheight{1.25}\smash{\begin{tabular}[t]{l}$\theta$\end{tabular}}}}%
    \put(0,0){\includegraphics[width=\unitlength,page=2]{AngleRhombus.pdf}}%
  \end{picture}%
\endgroup%
}
\end{center}
\caption{An isoradial domain and its various representations. To the left, we see the graph in red, in the middle, its dual graph superimposed in blue and the associated rhombic tiling in black and to the right, we see the rhombic tiling, the angle associated to an edge and two tracks of the graph. }
\label{fig:IsoReps}
\end{figure}

In this section, we introduce the concepts related to isoradial graphs which we will need throughout the paper. First, let us establish a bit of standard notation:~$\mathbf{Z}^2$ will henceforth always refer to the square lattice in its usual isoradial embedding, which is the rotation of~$\sqrt{2}\mathbb{Z}^2$ by~$\frac{\pi}{4}$. Set~$\Lambda_N$ to the Euclidean box~$[-N,N]^2$; when working with a fixed graph~$\bbG$, we will sometimes identify~$\La_N$ to the subgraph of~$\bbG$ it contains. 

One may note that any isoradial embedding of~$\mathbb{G}$ also gives rise to an isoradial embedding of its dual graph~$\mathbb{G}^*=(\mathbb{V}^*,\mathbb{E}^*)$ by embedding dual vertices as the centres of the circles inscribing the faces of~$\mathbb{G}$.

To any isoradial graph~$\mathbb{G},$ we associate a rhombic tiling of the plane~$\mathbb{G}^{\diamond}=(\mathbb{V}^{\diamond}
,\mathbb{E}^{\diamond})$ as follows:~$\mathbb{V}^{\diamond}=\mathbb{V}\cup \mathbb{V}^*$ and we draw an edge between every centre~$v^*\in \mathbb{V}^*$ of a face and each of the vertices~$v\in \mathbb{V}$ incident to that face. It is easy to check that~$\mathbb{G}^{\diamond}$ is an embedded planar graph and that each of its faces is a rhombus of side-length~$1$. Note that~$(\mathbf{Z}^2)^{\diamond}=\mathbb{Z}^2$.

Thus, to every isoradial graph, we can attach a certain rhombic tiling of the plane and any rhombic tiling gives rise to exactly two different isoradial graphs with the only ambiguity being which vertices to consider primal and which to consider dual. See Figure~\ref{fig:IsoReps}.

In this way, each edge~$e\in \mathbb{E}$ is contained in a unique rhombus~$R(e)\in (\mathbb{V}^{\diamond})^*$ and we define the angle~$\theta_e$ to be the interior angle of~$R(e)$ at one of the vertices opposite of~$e$ (which is well-defined by symmetry).

For a given rhombic tiling~$\mathbb{G}^{\diamond},$  define a \textbf{track} of~$\mathbb{G}^{\diamond}$ to be a sequence of parallel edges~$(e_n)_{n \in \mathbb{Z}}\subseteq \mathbb{E}^{\diamond}$ such that for every~$n$,~$e_n$ and~$e_{n+1}$ are edges of the same rhombus~$R\in (\mathbb{V}^{\diamond})^*$ while~$e_n$ and~$e_{n+2}$ are not. See Figure~\ref{fig:IsoReps}.
The transverse angle of the track~$(e_n)_{n \in \mathbb{Z}}$ is the angle between each of the~$e_n$ and the horizontal axis; we will generally take it to be in~$[0,\pi)$. 

We say that two distinct tracks~$t=(e_n)_{n\in \mathbb{N}}$ and~$t'=(f_n)_{n\in \mathbb{N}}$ \textbf{intersect} if there exist~$n,m$ such that~$e_n, f_m$ are part of the same rhombus~$R$. We say that two tracks~$t$ and~$t'$ are \textbf{parallel} if they do not intersect.

We will often identify a track with the mid-points~$z_n$ of the corresponding edges~$e_n$. We will also, without comment, refer to the tracks of~$\mathbb{G}^{\diamond}$ as the tracks of~$\mathbb{G}$ (and thus,~$\mathbb{G}$ and~$\mathbb{G}^*$ have the same tracks).

We say that a track~$t=(z_n)_{n\in \mathbb{Z}}$ has an \textbf{asymptotic direction} if there exists a vector~$v_t$ and a~$C>0$ such that, uniformly in~$n$ and~$m$,
$$
|z_n-z_m-v_t(n-m)| \leq C
$$

\begin{prop} \label{prop:asymp_direction}
Suppose that~$\mathbb{G}$ is bi-periodic. Then, every track of~$\mathbb{G}$ has an asymptotic direction.
\end{prop}

\begin{remark}\label{rem:asymp_direction}
It is worth noting that two tracks with non-parallel asymptotic directions necessarily intersect. As two tracks can intersect at most once, for bi-periodic graphs, tracks with parallel asymptotic directions cannot intersect. Therefore, the above proposition says that for any bi-periodic isoradial graph~$\mathbb{G}$, the tracks of~$\mathbb{G}$ can be partitioned into a finite number of track types labelled by their asymptotic directions.
\end{remark}

\begin{proof}
Let~$t=(z_n)_{n\in \mathbb{Z}}$ be a track of~$\mathbb{G}$ and let~$v$ and~$w$ be the vectors of translation corresponding to the actions of~$(1,0)$ and~$(0,1)$ respectively. By the pigeonhole principle, there exist distinct~$n,m\in \mathbb{Z}$ such that~$z_n$ and~$z_m$ represent the same class in~$\mathbb{G}/(\mathbb{Z}v\oplus \mathbb{Z} w)$. Accordingly, we have~$(j,k)\in \mathbb{Z}^2\setminus\{0\}$ such that
$z_n-z_m=jv+kw.$ It follows that~$t$ is invariant under translation by~$jv+kw$. It is straightforward to check that~$\frac{jv+kw}{n-m}$ is an asymptotic direction of~$t$.
\end{proof}

Finally, note that any biperiodic graph also has the bounded angles property \cite{ManoGrimmett} --- since the set of transverse angles of its tracks is finite --- and the square grid property \cite{ManoGrimmett}.

\subsection{The Star-Triangle Transform and the Track-Exchange Operator}

\begin{figure}
\begin{center}
\includegraphics[scale=1.2]{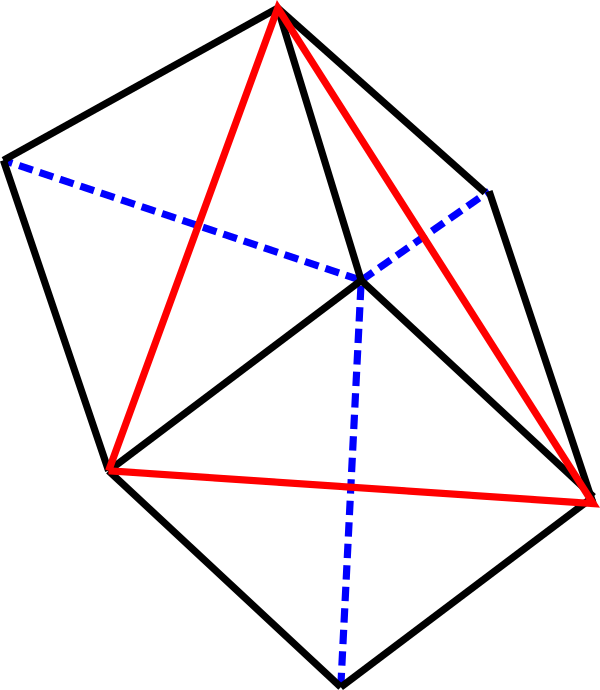} \hspace{2 cm} \includegraphics[scale=1.2]{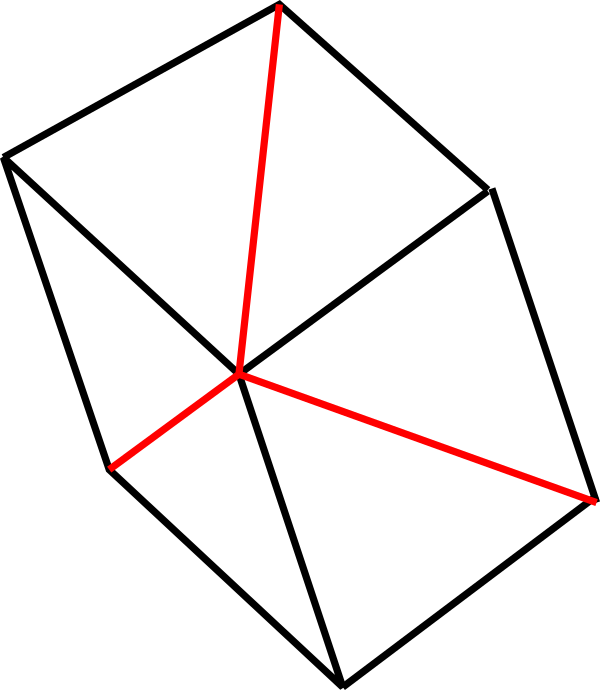}
\end{center}
\caption{To the left we see a triangle in red, its dual star in blue and the associated hexagon of rhombi in~$\mathbb{G}^{\diamond}$. To the right, the hexagon has been rotated by~$\pi$, replacing the triangle by a star. \label{Starntriangle}}
\end{figure}

The core ingredient in our treatment of universality among random-cluster models is the star-triangle transform (also known as the Yang-Baxter relation), which is a bit of integrability found in a staggering number of places throughout the planar statistical mechanics literature. It was first discovered by Kennelly in 1899 in the context of electrical networks \cite{electric} before appearing in Onsager's famous exact solution of the 2D Ising model \cite{Onsager}.

 In the present section, we shall focus on the purely graphical aspect of the star-triangle transformation, as that is going to guide the surgery procedures of the paper.

We call an isomorphic copy of the complete graph on three vertices a \textbf{triangle}. The set of dual edges to a triangle is called a \textbf{star}. One might note that, in the rhombic tiling description of an isoradial graph, this corresponds to three neighbouring rhombi forming a hexagon. See Figure~\ref{Starntriangle}.

For any isoradial graph~$\mathbb{G}$ and any star or triangle in~$\mathbb{G}$, we can define the star-triangle transform~$\mathbf{S}(\mathbb{G})$, which replaces the star by a triangle or vice versa. From the rhombic tiling point of view, this corresponds to rotating the three rhombi forming the hexagon by~$\pi$. Thus, the star-triangle transform preserves isoradiality of the graph.

\begin{figure} 
\begin{center}
\includegraphics[width = .25\textwidth, page = 1]{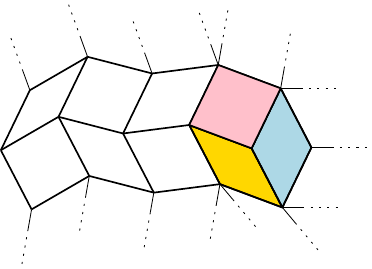} \hspace{0.05\textwidth}
\includegraphics[width = .25\textwidth, page = 2]{track_exchange} \hspace{0.05\textwidth}
\includegraphics[width = .25\textwidth, page = 3]{track_exchange} 
\end{center}
\caption{A finite track-exchange in a wider isoradial graph $\bbG$. Successive star-triangle transforms applied in the coloured hexagon move the blue rhombus to the left, thus exchanging the two horizontal rows of rhombi.}
\label{fig:TrackExchange}
\end{figure}

Now, if two tracks~$t$ and~$s$ run next to each other before intersecting, we may apply multiple consecutive star triangle transforms to locally change the relative position of the two tracks. Such a sequence of transformations is called a finite \textbf{track-exchange}; see  Figure~\ref{fig:TrackExchange}. As in \cite{Rotation}, infinite track-exchanges may also be defined; see also Remark~\ref{rem:infinite_track_ex} below. 

\subsection{The Random-Cluster Model}

We now introduce the random-cluster model on isoradial graphs in more detail. We also recall its elementary properties that we are going to need later on. We refer to the monographs \cite{DC17,Gri06} for thorough introductions to the random-cluster model.

\begin{definition} For a finite subgraph~$G=(V,E)$ of an isoradial graph~$\mathbb{G}$,~$q\geq 1$, vector~$\boldsymbol{p}=p_e\in [0,1]^E$ and partition~$\xi$ of the vertex boundary of~$G$, the \textbf{random-cluster measure} on~$G$ with boundary condition~$\xi$ is the measure~$\phi^{\xi}_{\boldsymbol{p},q,G}$ on~$\{0,1\}^E$ given by
$$
\phi^{\xi}_{\boldsymbol{p},q,G}[\{\omega\}]=\frac{1}{Z_{\boldsymbol{p},q,G}^{\xi}}q^{\kappa^{\xi}(\omega)}\prod_{e\in E} p_e^{\omega(e)} (1-p_e)^{1-\omega(e)},
$$
where~$\kappa^{\xi}$ is the number of connected components in the graph~$$
G^{\xi}_{\omega}=(V, E_{\omega})/\sim_{\xi},
$$
~$E_{\omega}:=\{ e\in E|\;\omega(e)=1\},$~$\sim_{\xi}$ is the equivalence relation with equivalence classes given by~$\xi$ and~$Z_{\boldsymbol{p},q,G}^{\xi}$ is a normalising constant called the partition function.
 
~$E_{\omega}$ is called the set of \textbf{open edges} and~$E\setminus E_{\omega}$ is called the set of \textbf{closed edges}.
\end{definition}
\begin{prop}[Domain Markov Property] \label{DMP}
If~$G_1=(V_1,E_1)\subseteq G_2=(V_2,E_2)$ are two finite subgraphs of an isoradial graph~$\mathbb{G}$, we write~$\omega_1:=\omega|_{E_1}$ and~$\omega_2:=\omega|_{E_2\setminus E_1}$ and have
$$
\phi_{\boldsymbol{p},q,G_2}^{\xi}[\omega_1\in A|\;\omega_2]=\phi_{\boldsymbol{p},q,G_1}^{\xi_{\omega_2}}[A],
$$
where~$v$ and~$w$ belong to the same class of~$\xi_{\omega_2}$ if and only if they are connected by a path in~$((V_2\setminus V_1), E_{\omega_2})/\sim_\xi$.
\end{prop}
\begin{definition}
For an isoradial graph~$\mathbb{G}=(\mathbb{V},\mathbb{E})$, we say that a probability measure on~$\phi_{\boldsymbol{p},q}$ on~$\{0,1\}^{\mathbb{E}}$ is an \textbf{infinite-volume random-cluster measure} on~$\mathbb{G}$ if, for any finite subgraph~$G=(V,E)$ of~$\mathbb{G}$, we write~$\omega_1=\omega|_E$ and~$\omega_2=\omega|_{E^c}$ and almost surely have
$$
\phi_{\boldsymbol{p},q}[\omega_1\in A|\; \omega_2]=\phi_{\boldsymbol{p},q,G}^{\xi_{\omega_2}}[A],
$$
where~$v$ and~$w$ belong to the same element of~$\xi_{\omega_2}$ if and only if they are connected by a path in~$(\mathbb{V}\setminus V,E_{\omega_2})$.
\end{definition}

For an isoradial graph~$\bbG = (\bbV,\bbE)$ and~$\omega \in \{0,1\}^\bbE$ a FK-percolation configuration. 
We define the {\em dual configuration}~$\omega^*$ on~$\bbG^*$ (defined in Section~\ref{sec:isoradial_graphs}) by 
\begin{align}
\omega^*(e^*) = 1 - \omega(e) \qquad \forall e \in \bbE,
\end{align}
were~$e^*$ is the edge dual to~$e$. 

\begin{figure}
\begin{center}
\includegraphics[width = .4\textwidth]{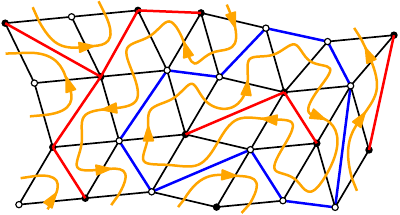}
\caption{An isoradial graph ($\bbG^\diamond$ is in black), with a primal configuration (in red), its dual configuration (in blue) and the loop representation (in orange).}
\label{fig:loop_rep}
\end{center}
\end{figure}

The medial graph of~$\bbG$ is the dual of~$\bbG^\diamond$. 
Its vertices are the intersection points between the primal and dual edges of~$\bbG$; vertices of the medial graph are connected if they are centres of adjacent rhombi.
The loop representation of~$\omega$ is the partition of the edges of the medial graph into loops and bi-infinite paths
in such a way that no path or loop intersects any open edge of~$\omega$ or open dual edge of~$\omega^*$. 
These may be represented as in Figure~\ref{fig:loop_rep} using bent edges and 
small turns at the vertices of the medial graph. 

The loops are oriented to have the dual open edges on their left. 
Throughout this paper, the loop representation of~$\omega$ contains no infinite paths; it is formed entirely of finite loops. 

\subsection{Isoradial weights and the Star-Triangle Transform}

 In this paper, we shall only be interested in so-called critical models with isoradial weights. These weights are chosen so as to be compatible with the isoradial embedding of the graph and the star-triangle transform. For a given star-triangle transform~$\mathbf{S}$ and~$\omega\in \{0,1\}^E$, we can associate a random configuration~$\mathbf{S}(\omega)\in \{0,1\}^{\mathbf{S}(E)}$ as follows:
\begin{enumerate}
\item[$(i)$]~$\omega(e)=\mathbf{S}(\omega)(e)$ for every~$e$ not affected by~$\mathbf{S}$.
\item[$(ii)$] For the edges affected by~$\mathbf{S},$ we pair configurations as indicated in Figure~\ref{Pairing}. In cases of ambiguity, we select the result with probabilities proportional to the random-cluster weights\footnote{One may note that the relative number of components in each case of ambiguity is independent of the boundary condition}. Solid red lines are open edges and dashed blue ones are closed edges. \\
\begin{figure}
\begin{center}
\includegraphics{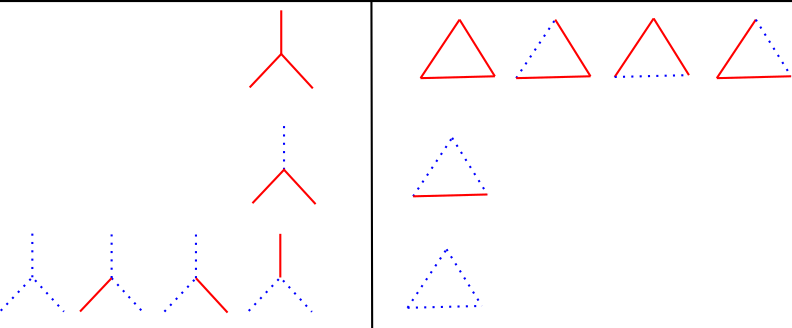}
\end{center}
\caption{  \label{Pairing}The pairing of percolation configurations under the star-triangle transform. We have omitted symmetric cases.}
\end{figure}

\end{enumerate}
\begin{remark}
We may note that the connection field of~$\omega$ is preserved, meaning that for~$v,w\in \mathbb{V}\cap \mathbf{S}(\mathbb{V})$,~$v$ is connected to~$w$ in~$\omega$ if and only if they are in~$\mathbf{S}(\omega)$.
\end{remark}
Define the isoradial weights
$$
f^q(\theta)=\begin{cases}
\frac{\sqrt{q}\sin(r(\pi-\theta))}{\sin(r\theta)+\sqrt{q}\sin(r(\phi-\theta))} & q<4\\
\frac{2\pi-2\theta}{2\pi-\theta} & q=4
\end{cases},
$$
where~$r=\frac{1}{\pi}\cos^{-1}(\sqrt{q}/2)$.
The following is essentially due to Kenyon \cite{Kenyon}, but the precise statement is to be found in \cite{JiHuangLi}.

\begin{prop}[\cite{JiHuangLi}, Proposition 2.3]
For every~$q\in [1,4]$, any finite part~$G$ of an isoradial graph~$\mathbb{G}$ and any boundary condition~$\xi$, the star-triangle transform preserves the random-cluster measure with parameters~$p^{\mathbb{G}}_e=f^q(\theta_e)$. I.e., for these edge weights, we have
$$
\mathbf{S}(\phi_{{G}}^{\xi})=\phi^{\xi}_{\mathbf{S}({G})}.
$$
\end{prop}

\begin{remark} 
The star-triangle transform neither removes nor adds new boundary vertices so that it does, indeed, make sense to regard the same~$\xi$ on both sides of the equation.
\end{remark}

Throughout the rest of the paper, we will assume that~$q\in [1, 4]$ is fixed. In this range of parameters, it was shown in \cite{JiHuangLi} that for bi-periodic isoradial graphs~$\bbG$, there is a unique infinite-volume measure, which we denote by~$\phi_{\bbG}$.  
When~$\bbG$ is not bi-periodic, we will require that it has a unique infinite volume measure --- we expect this to be the case for all isoradial graphs in this regime. 

\begin{remark}\label{rem:infinite_track_ex}
As shown in \cite{Rotation}, if an infinite isoradial graph~$\bbG$ has a unique infinite volume measure~$\phi_\bbG$ and contains two parallel tracks~$s$ and~$t$ that run next to each other (i.e. no rhombus is contained strictly between~$s$ and~$t$), one may perform an {\bf infinite track-exchange} of these two tracks. 

More precisely, if we write~$\bbG'$ for the isoradial graph obtained by exchanging the rhombi of~$t$ and~$s$, there exists a random map~${\mathbf T}$ between configurations on~$\bbG$ and configurations on~$\bbG'$ such that 
\begin{align}
\omega \sim \phi_{\bbG} \quad \Longrightarrow\quad {\mathbf T}( \omega) \sim \phi_{\bbG'},
\end{align}
where~$\phi_{\bbG'}$ is the unique infinite volume measure on~$\bbG'$. 
Furthermore,~$\mathbf T$ preserves the connectivity between vertices common to~$\bbG$ and~$\bbG'$ -- that is, all vertices except those contained between~$s$ and~$t$. 
Indeed,~${\mathbf T}$ may be viewed as a bi-infinite sequence of successive star-triangle transformation. Direct constructions of~${\mathbf T}$ are also possible, see for instance~\cite{PelRom15}. 
\end{remark}

Let us briefly remind the reader of the RSW theory, which has played a central role in the study of planar critical percolation models.  

\begin{definition}\label{def:RSW}
Let~$Q$ be rotation and translation of~$[0,4n]\times [0,n]$ for~$n\geq 1$. Let~$\mathscr{C}(Q)$ be the event that there exists a path of open edges crossing~$Q$ between its two short sides.
We say that an isoradial graph~$\mathbb{G}$ has the \textbf{strong RSW property} with parameter~$c>0$ if, 
for every~$Q$ as above and~$n\in \mathbb{N}$, every event~$B$ depending on the edges at distance at least~$n$ from~$Q$, we have  
\begin{align}\label{eq:RSW}
	\phi_{\mathbb{G}}[\mathscr{C}(Q)\,|\, B]\in [c,1-c]
\end{align}
We write~$\mathbb{G}\in \mathcal{RSW}(c)$ if the above holds for~$\bbG$. 
\end{definition}

The results in \cite{JiHuangLi} imply that for any bi-periodic isoradial graph~$\bbG$,~$\phi_{\mathbb{G}}$ admits a strong RSW property. One standard consequence of RSW is a mixing property. See e.g. \cite[Proposition 2.9]{ScalingRelations}.
\begin{prop}[Mixing] \label{mixing}
For every~$c>0$, there exists~$C>0$ such that for every~$\mathbb{G}\in \mathcal{RSW}(c)$ and events~$A$ depending only on edges in~$\Lambda_r$ and~$B$ depending only on edges outside of~$\Lambda_R$, we have
$$
\big|\phi_{\mathbb{G}}[A\cap B]-\phi_{\mathbb{G}}[A]\phi_{\mathbb{G}}[B]\big|
\leq C \left(\tfrac{r}{R}\right)^C \phi_{\mathbb{G}}[A]\phi_{\mathbb{G}}[B].
$$
\end{prop}

\subsection{Extensions of isoradial domains} 

An isoradial domain~$\calD$ is a finite, simply connected set of rhombi of a rhombic tiling~$\bbG^\diamond$. 
We call the isoradial graph~$\mathbb{G}$ an {\em extension} of~$\calD$. 
It is called a {\em natural extension} if it is such that there exists a unique infinite volume measure~$\phi_\bbG$ on~$\mathbb G$. In particular, biperiodic graphs are natural extensions.

\begin{prop}\label{prop:natural_extension}
	Let~$\calD$ be an isoradial domain and~$\bbG$ and~$\bbG'$ be two natural extensions of~$\calD$. 
	Then, the marginals of~$\phi_\bbG$ and~$\phi_{\bbG'}$ on~$\calD$ are the same. 
	In fact, even the distribution of the boundary conditions on~$\calD$ induced by the outside configuration is the same under both measures.
\end{prop}

\begin{proof}
	The proof has two parts. 
	We first define a particular extension~$\bbG_0$ of~$\calD$; we will construct it by adding rhombi to~$\calD$ sequentially. 
	There are two stages to the construction: first we {\em convexify}~$\calD$, and then we complete it into an isoradial graph in the full plane. 

	Write~$t_1,\dots, t_k$ for the tracks traversing~$\calD$. Recall that any rhombus of an isoradial graph corresponds to the intersection of two tracks. Thus, the potential rhombi formed by the tracks~$t_1,\dots, t_k$ is that of the pairs~$\{t_i,t_j\}$ with~$i\neq j$. Some of these intersections may appear in~$\calD$, but not necessarily all of them. If~$\{t_i,t_j\}$ corresponds to no rhombus of~$\calD$, we say that the intersection is {\em missing}. We say that~$\{t_i,t_j\}$ is {\em accessible} if it corresponds to a rhombus which shares two edges with existing adjacent rhombi on the boundary of~$\calD$.

	\begin{figure}
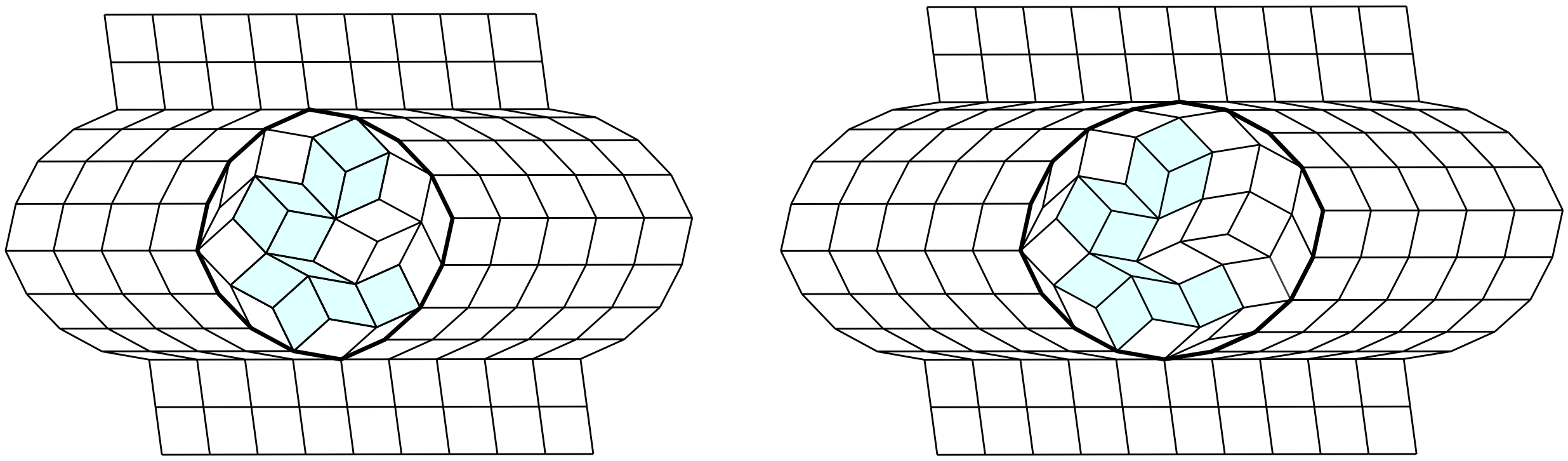

	\begin{center}
	\includegraphics[width = 0.48\textwidth, page = 2]{completion.pdf}\quad
	\includegraphics[width = 0.48\textwidth, page = 4]{completion.pdf}
	\caption{{\em Left:} The default extension~$\bbG_0$ of an isoradial domain~$\calD$. 
	The convexification~$\calD'$ is surrounded by the bold line. 
	{\em Right:} A different extension~$\bbG_\Lambda$ of~$\calD$. The convexified part~$\Lambda'$ is also highlighted; outside of it we use the same track angles as for~$\bbG_0$. }
	\label{fig:completion1}
	\end{center}
	\end{figure}
	
	Consider the lexicographical order on the pairs~$\{t_i,t_j\}$ with~$i <  j$. 
	Let~$\{t_i,t_j\}$ denote the smallest pair so that~$\{t_i,t_j\}$ is both missing and accessible. If such a pair exists, add the rhombus and restart the procedure. By repeatedly doing this, we arrive at a stage where no such pair of tracks exists. Let~$\calD'$ denote the isoradial domain obtained by this procedure. 
	
	We claim that the union of the rhombi of~$\calD'$ form a convex set, and that~$\calD'$ contains all intersections between tracks of distinct transverse angles.
	Indeed, if no rhombus as above may be added, the edges forming the boundary of~$\calD'$ have angles in~$[0,2\pi)$ that form an increasing sequence along the boundary in counterclockwise order (if starting from the lowest angle). In other words, the boundary can be parametrised by a curve which only ever turns to the left, and so,~$\calD'$ is convex. Moreover, it implies that any track enters and exits~$\calD'$ at opposing sides. For this to be true for all finitely many angles appearing, each angle must intersect every other angle. See Figure~\ref{fig:completion1} for an illustration. 
	
	Finally, we construct~$\bbG_0$. Consider two angles~$\theta_1\neq\theta_2$ 
	so that none of the tracks~$t_1,\dots, t_N$ has these particular transverse angles.  
	Complete~$\calD'$ with tracks of transverse angles~$\theta_1$ and~$\theta_2$ so that the tracks~$t_1,\dots,t_k$ intersect all those with angle~$\theta_1$, but none of those with angle~$\theta_2$. 
	
	For simplicity of exposition, we will consider~$\theta_1 = 0$ henceforth. 
	We call the tracks of transverse angle~$\theta_1$ {\em vertical} and those of transverse angle~$\theta_2$ {\em horizontal}.
	
	Notice that~$\bbG_0$ has the structure of a square lattice outside of~$\calD'$. As a result, the RSW estimates of \cite[Thm. 1.1]{JiHuangLi} apply for rectangles sufficiently far from~$\calD'$, which implies that~$\bbG_0$ has a unique infinite volume FK-percolation measure. 
	Thus,~$\bbG_0$ is a natural extension of~$\calD$.  We will consider~$\bbG_0$ to be the {\em default} extension of~$\bbD$. 
	\bigskip 
	
	In the second part of the proof, we show that any other natural extension of~$\calD$ produces the same measure on~$\calD$ as~$\phi_{\bbG_0}$. 
	Fix such a natural extension~$\bbG_1$. 
	Let~$\omega_0$ be a configuration on~$\calD$ and~$\eps > 0$. 
	Let~$\Lambda$ be a large enough connected part of~$\bbG_1$ so that 
	\begin{align}\label{eq:lala0101}
		\max_{\xi,\xi'}\big|\phi_\Lambda^{\xi}[\omega = \omega_0 \text{ on~$\calD$}] - \phi_\Lambda^{\xi'}[\omega = \omega_0 \text{ on~$\calD$}]\big| \leq \eps.
	\end{align}
	Define now~$\bbG_\Lambda$ as the completion of~$\Lambda$ obtained as above, by first convexifying~$\Lambda$ to obtain~$\Lambda'$, then completing it with tracks of transverse angles~$\theta_1$ and~$\theta_2$. In light of \eqref{eq:lala0101}, we find 
	\begin{align}\label{eq:lala0102}
		|\phi_{\bbG_1}[\omega = \omega_0 \text{ on~$\calD$}] - \phi_{\bbG_\Lambda}[\omega = \omega_0 \text{ on~$\calD$}]| \leq \eps.
	\end{align}

	Next, we apply star-triangle transformations to~$\bbG_\Lambda$ so as to create an arbitrarily large portion around~$\calD$ identical to~$\bbG_0$. These transformations will not affect~$\calD$. 
	This is done in two steps. \medskip 
	
	\noindent {\bf Step 1:} We first show that by applying star-triangle transformations to~$\bbG_\Lambda$ outside of~$\calD$,
	we may reconstruct~$\calD'$. 
	See Figure~\ref{fig:completion2} for an illustration of this step.

	\begin{figure}
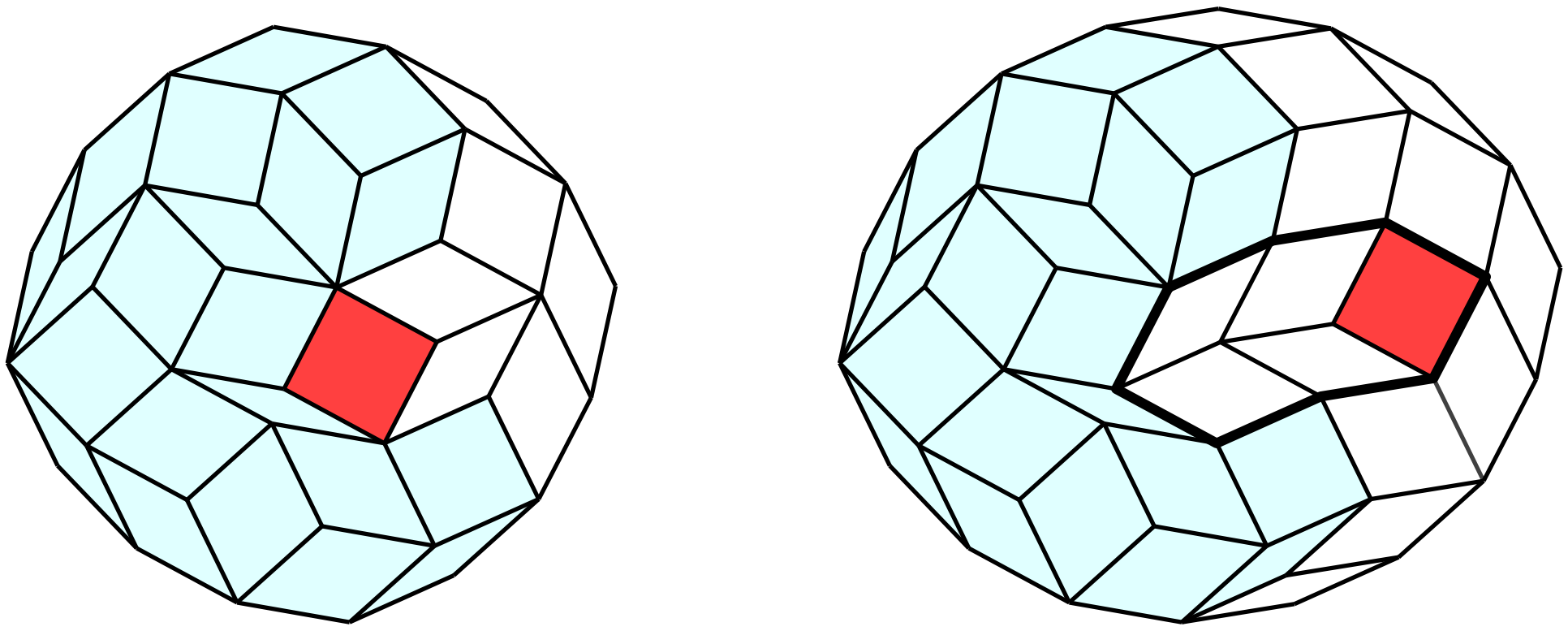

	\begin{center}
	\includegraphics[width = 0.35\textwidth, page = 5]{completion2.pdf}\qquad
	\includegraphics[width = 0.35\textwidth, page = 6]{completion2.pdf}
	\caption{The procedure used in Step 1 of the transformation of~$\bbG_\Lambda$ into~$\bbG_0$. 
	{\em Left:} the part~$\calD'$ of~$\bbG_0$ with~$\calA$ in blue and the rhombus~$R_j$ marked in red. 
	{\em Right:} the convex part~$\Lambda'$ of~$\bbG_\Lambda$ with the rhombus~$\tilde R$ in red. 
	This is the vertex which needs to be moved to the position occupied by~$R$; we do so by applying star triangle transformations in the delimited region. Here, two transformations suffice.}
	\label{fig:completion2}
	\end{center}
	\end{figure}

	Write~$R_1,R_2,\dots$ for the rhombi of~$\calD'\setminus \calD$ in the order in which they were added when convexifying~$\calD$. 
	Let~$j$ be maximal so that~$R_1, \dots, R_{j-1}$ appear in~$\Lambda'$ at the same position as in~$\calD'$. 
	Finally, write~$\calA = \calD \cup R_1\cup \dots \cup R_{j-1}$. 
	
	If~$\calA = \calD',$ we are done, and we may move on to Step 2. 
	Otherwise, observe that~$R_j$ is the intersection of two tracks,~$s$ and~$t,$ which cross~$\calD$. 
	These two tracks necessarily intersect in~$\Lambda'$ (as~$\Lambda'$ is convex), and we denote by~$\tilde R$ the rhombus of~$\bbG_\Lambda$ corresponding to their intersection. 

	Since~$\calA$ is connected, there exists a finite region in~$\bbG_\Lambda$  delimited by~$s$,~$t$ and~$\calA$. Write~$\calB$ for the set of rhombi contained in this region, including those contained on~$s$ and~$t$, but not those contained in~$\calA$.
	
	Note that since~$R_j \neq \tilde R$ by assumption,~$|\mathcal{B}|>1$. 
	Therefore, there necessarily exists a hexagon that, when flipped, reduces the number of rhombi in~$\calB$ by at least one.
	This transformation involves the rhombi in~$\calB$, and therefore does not affect~$\calA$, and in particular not~$\calD$.  

	Applying such transformations repeatedly, we reduce~$\calB$ to a single rhombus, which in turn yields a new graph~$\tilde{\mathbb{G}}_{\Lambda}$ with ~$\tilde{\calA}=\calA\cup \{R_j\}$.
	By repeated applications of this argument,~$\calA$ may be increased until it is equal to~$\calD'$.\medskip 

	\noindent {\bf Step 2:}	We may now assume that~$\bbG_\Lambda$ contains~$\calD'$, and we wish to apply star-triangle transformations to it, outside of~$\calD'$, so as to recreate an arbitrarily large portion of~$\bbG_0$.
	See Figure~\ref{fig:completion3} for an illustration of this step. 
	
	For ease of exposition, we will assume that we may perform infinitely many star-triangle transformations  and infinite track-exchanges, 
	and will transform~$\bbG_\La$ into~$\bbG_0$ exactly.
	This may be rendered formal and re-interpreted to show how~$\bbG_\La$ may be transformed into~$\bbG_0$ on an arbitrarily large region by a finite number of transformations. 
	
	\begin{figure}
	\begin{center}
	\includegraphics[width = 0.65\textwidth, page = 7]{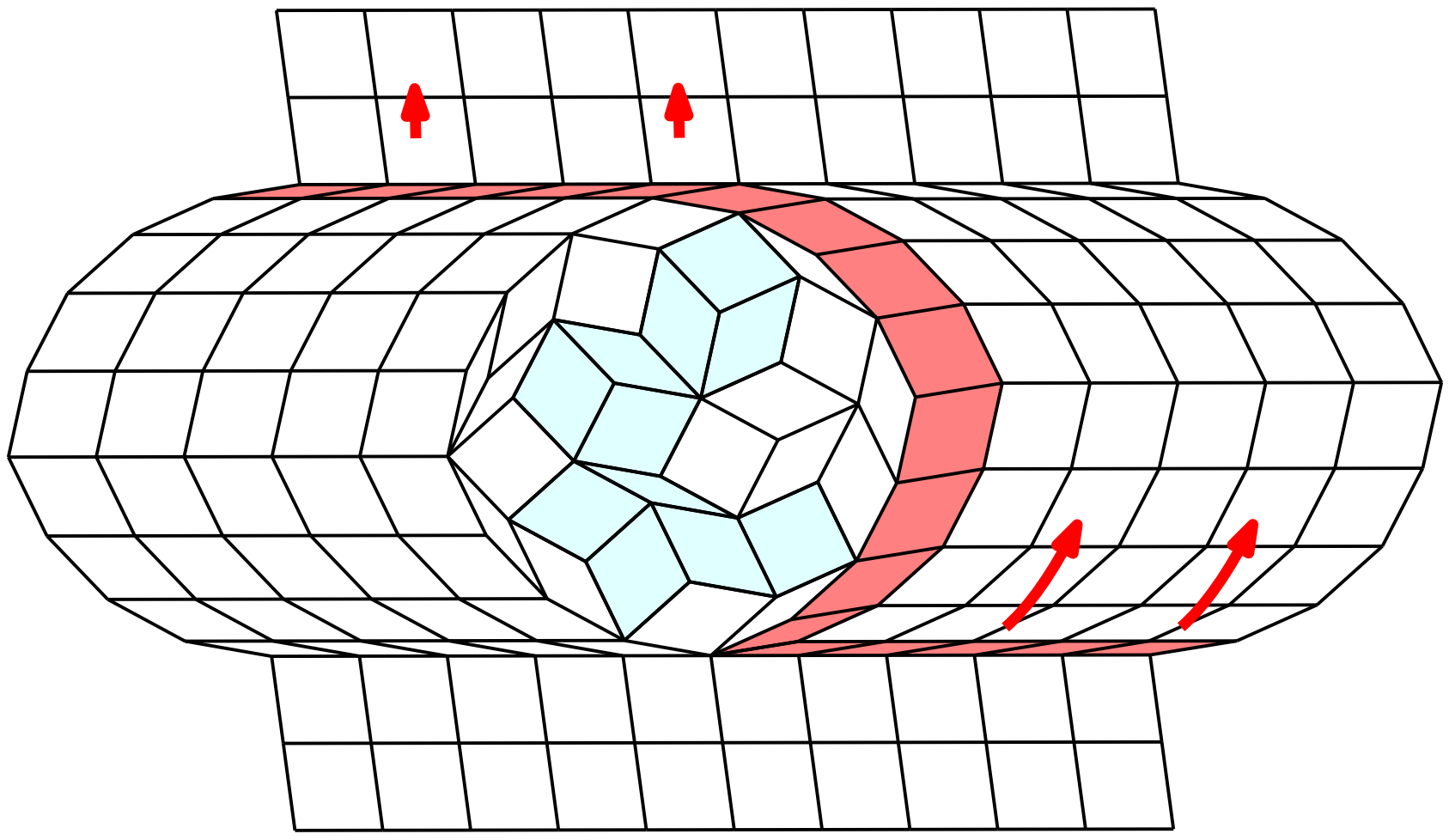}
	\caption{The procedure used in Step 2 of the transformation of~$\bbG_\Lambda$ into~$\bbG_0$. 
	The graph~$\bbG_\Lambda$ may now be considered to contain~$\calD'$ within the convex region~$\Lambda'$. 
	The tracks that intersect~$\Lambda'$ but not~$\calD'$ (here, there is only one, marked in red) need to be sent to infinity. 
	We do this by star-triangle transformation and track-exchanges. Here the red track may be moved upwards indefinitely.}
	\label{fig:completion3}
	\end{center}
	\end{figure}
	
	We refer to the tracks of~$\bbG_\La$ that cross~$\Lambda'$ but not~$\calD'$ as {\em extra tracks}. 
	If no such tracks exist, then~$\calD' = \Lambda'$, and~$\bbG_\Lambda=\bbG_0$. 
	Assume next that at least one extra track exists. We will eliminate these tracks by ``pushing them to infinity'' using star-triangle transformations  and track-exchanges.
	
	As extra tracks cross neither~$\calD'$ nor any of the horizontal tracks, they either run above or below~$\calD'$. 
	We will limit ourselves to eliminating those running above~$\calD'$. 
	Among these, consider the extra track~$t$ with minimal transverse angle; 
	this is the extra track which is topmost to the left of~$\La'$ and bottommost to the right.
	We write~$s_1,s_2,\dots$ for the horizontal tracks above~$\La'$, in vertically increasing order.

	Applying star-triangle transformations on the right of~$t$, this track may be moved to run adjacent to~$s_1$ by a procedure which is very similar to the one used in step 1.
	Then, applying track-exchanges between (the image of)~$t$ and the horizontal tracks~$s_1,s_2\dots$,~$t$ may be sent upwards to infinity.
	This concludes Step 2.  
	\bigskip

	Observe now that none of the star-triangle transformations and track-exchanges applied above to~$\bbG_\Lambda$ affects the marginal of the measure in~$\calD$ or even the boundary conditions on that domain (i.e. the connections between points on the boundary). Moreover, since these transformations modify ~$\bbG_\Lambda$ so as to contain an arbitrary large region around~$\calD$ that is identical to~$\bbG_0$, and since the latter has a unique infinite-volume measure, we conclude that 
	\begin{align}\label{eq:lala0103}
		\phi_{\bbG_\Lambda}[\omega = \omega_0 \text{ on~$\calD$}] = \phi_{\bbG_0}[\omega = \omega_0 \text{ on~$\calD$}].
	\end{align}
	Combining the above with \eqref{eq:lala0102} and using the fact that~$\eps$ is arbitrary, we conclude that 
	\begin{align}\label{eq:lala0104}
		\phi_{\bbG_1}[\omega = \omega_0 \text{ on~$\calD$}] = \phi_{\bbG_0}[\omega = \omega_0 \text{ on~$\calD$}],
	\end{align}
	as desired. The same applies to the boundary conditions induced on~$\calD$ by the configuration outside.
\end{proof}

\subsection{Irreducibility of pavings of convex domains}

We close this section with a useful combinatorial observation. 

\begin{lemma}\label{lem:irred_convex}
Let $\calD$ and $\calD'$ be two isoradial domains that occupy the same convex set of $\bbR^2$. 
Write $2M$ for the length of the boundary of $\calD$ and $\calD'$. 
There exists a sequence $\bfS_1,\dots, \bfS_L$ of at most $M^3$ star triangle transformations that may be applied in order to $\calD$ and such that 
\begin{align}\label{eq:irred_convex}
(\bfS_L \circ \cdots \circ \bfS_1) (\calD)= \calD'.
\end{align}
\end{lemma}

\begin{proof}
	It suffices to prove \eqref{eq:irred_convex} for a specific choice of $\calD'$, with $L \leq M^3/2$. 
	This specific choice of $\calD'$ is defined as follows:
	
	Write $\calC$ for the convex set occupied by $\calD$
	and $e_1,\dots, e_{2M}$ for the edges of the boundary of $\calC$, in counterclockwise order, so that $e_1,\dots, e_M$ have increasing angles $\theta_1,\dots \theta_M$ in $[0,\pi)$.
	For any paving of $\calC$ with rhombi, write $t_i$ for the track entering through $e_i$. Such a track necessarily exits through $e_{M+k+i_0-(i-i_0)},$ where $i_0$ is the first index with $\theta_{i_0}=\theta_i$ and $k+i_0$ is the last index with $\theta_{k+i_0}=\theta_i$.
	
	We define the standard paving $\calD'$ of $\calC$ inductively. Set $i_0$ equal to smallest index $j$ such that $\theta_j=\theta_M$ and set $k=M-i_0$.
	We first require that $t_M$ runs along the left boundary of $\calC$ from $e_M$ to $e_{2M-k}$.
	This fixes $M-k-1$ rhombi in $\calC$. 
	The rest of $\calC$ is surrounded by a convex perimeter with angles $\theta_1,\dots, \theta_{M-1}, \theta_1, \dots, \theta_{M-1}$. 
	Proceeding by induction, this constructs the standard paving.
	
	Now that $\calD'$ is defined, we explain how to pass from $\calD$ to $\calD'$ using star-triangle transformations. 
	This is again done by induction. 
	Consider the track $t_M$ of $\calD$. If this is in the same position as in $\calD'$, fix this region and proceed with transforming the rest of the domain, which corresponds to the case $M-1$.
	
	If $t_M$ in $\calD$ differs from its position in $\calD'$, ``push'' it to the left. 
	Indeed, our assumption implies the existence of at least one rhombus in $\calD$ which lies strictly to the left of $t_M$.
	One may then find  such a rhombus $R$ that has two sides adjacent to $t_M$. For instance, one may orient the tracks $t_{1},\dots,t_{M-k}$ from the left of $t_M$ towards the right, inducing a partial order on rhombi to the left of $t_M$. Then, pick $R$ to be any maximal rhombus of this partial ordering. See Figure \ref{fig:Irreducible_Paving}.
	Applying a star triangle-transform to $R$ and its two neighbouring rhombi on $t_M$, we obtain a graph $\bfS(\calD)$ with $R$ to the right $t_M$. 
	Repeating this procedure pushes $t_M$ to the left, in the same position as in $\calD'$. 
	
	\begin{figure}
	\begin{center}
	\includegraphics{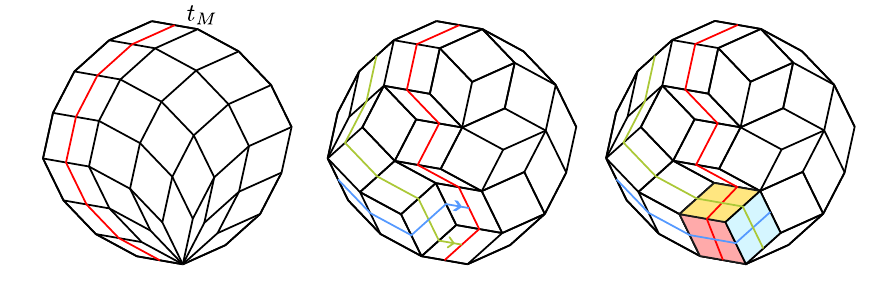}
	\label{fig:Irreducible_Paving}
	\end{center}
\vspace{-1cm}	\caption{A convex set $\mathcal{C}$ with its standard paving on the left. In the middle, another paving of $\mathcal{C}$, where $t_M$ is not the leftmost track. This implies the existence of a hexagon which may be flipped with the star-triangle transform and moved to the right of $t_M$, such as the intersection of the blue and green track, which is maximal for the ordering of rhombi given by the orientations of the tracks $t_1,\dots t_{M-1}$ from left to right.}
\end{figure}		
	
	The number of star-triangle transformations needed to push $t_M$ to the left is bounded by the overall number of rhombi in $\calD$,
	which is at most $M(M-1)/2$ - corresponding to the case where all angles are distinct. 
	Since this procedure is then applied to $t_{M-1}$ in the  remaining part, then to $t_{M-2}$ etc, we conclude that the total number of star-triangle transformations needed to pass from $\calD$ to $\calD'$ is bounded above by 
	\begin{align}
	\sum_{l=1}^M \frac{l(l-1)}2 \leq \frac{M^3}2.
	\end{align}
	This concludes the proof.
\end{proof}

\section{The Cube Flip Lemma}\label{sec:cube_flip}

In the context of isoradial graphs, the star-triangle transformation may be viewed as the permutation of three adjacent rhombi which together form a hexagon. 
Here, we consider a large-scale version of this transformation and describe its effect on the configuration. 

\begin{figure}
\begin{center}
\includegraphics{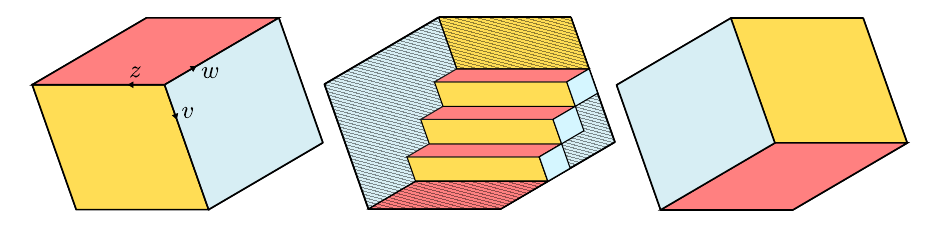}
\end{center}
\vspace{-1cm}
\caption{\label{CubePic} {\em Left:} a cube with the block of old rhombi in yellow, the new ones in red and the transverse ones in blue. 
{\em Right:} the flipped cube, with the three blocks exchanged.
{\em Middle:} a stage~$\bbG_{kN}$ in the transformation from the initial to the flipped cube. The frozen block is hashed and the active one appears bright. To get to~$\bbG_{(k+1)N}$, all new rhombi in the active block are pushed down by track-exchanges. The blue transverse rhombi are also pushed to the left side.
}
\end{figure}

	Given three unit vectors~$v,w$ and~$z$ not contained in the same closed half-plane, we define the {\em cube}~$\mathcal{C}(v,w,z,N)$ to be the isoradial domain obtained by gluing translates of three~$N\times N$ rhombi, each paved with~$N^2$ unit rhombi with sides~$v\times w$,~$w\times z$ and~$z\times v$, respectively, as in Figure~\ref{CubePic}, left diagram. 
	Furthermore, we define the \textbf{flipped cube}~$\mathcal{C}^{\mathfrak{F}}(v,w,z,N)=\mathcal{C}(-v,-w,-z,N)$.
If~$\bbG$ is an extension of~$\mathcal{C}(v,w,z,N)$, write~$\mathbb{G}^{\mathfrak{F}}$ for the same graph, but with~$\mathcal{C}(v,w,z,N)$ replaced by the flipped cube ~$\mathcal{C}^{\mathfrak{F}}(v,w,z,N)$.

As~$v,w,z$ and~$N$ will mostly be fixed in this section, we will omit them from notation when no ambiguity is possible. The purpose of this section is to prove the following:

\begin{thm}[Cube flip]\label{thm:cube_flip}  
	For any three unit vectors~$z,v,w$ not contained in the same closed half-plane, 
	there exists a constant~$c>0$ with the following property:
	For any~$N\geq1$ and any natural extension~$\mathbb{G}$ of~$\mathcal{C}(v,w,z,N)$, there exists a coupling~$\mathbb{P}$ between~$\omega_1\sim \phi_{\mathbb{G}}$ and~$\omega_2\sim \phi_{\mathbb{G}^{\mathfrak{F}}}$ such that
	\begin{itemize}
		\item[(i)]~$\omega_1(e) = \omega_2(e)$ for any edge~$e$ outside of~$\mathcal{C}$;
		\item[(ii)] the connections induced by~$\omega_1 \cap \mathcal{C}$ on the vertices of~$\partial \mathcal{C}$
		are the same as those induced by~$\omega_2 \cap \mathcal{C}^{\mathfrak{F}}$;
		\item[(iii)]~$\displaystyle \mathbb{P}\big[ {d}_{\rm global\, CN}\big({\omega_1},{\omega_2}\big)\geq N^{1-c}\big]\leq \tfrac{1}{N^c}$. 
	\end{itemize}
\end{thm}

The proof of the above requires a slight re-working of the results of \cite{Rotation}. Indeed, the flipping of the cube from~$\mathcal{C}$
to~$\mathcal{C}^{\mathfrak{F}}$ may be done by a series of star-triangle transformations that resemble the track-exchanges of \cite{Rotation}. 
Some additional details need to be considered due to the boundary effects that appear when flipping the cube, rather than infinite horizontal tracks. 
Given the complexity of the argument in \cite{Rotation}, we will give a short sketch below. 
The rest of the section is dedicated to the proof of Theorem~\ref{thm:cube_flip}.

\subsection{RSW theory in the cube}

In the proof of Theorem~\ref{thm:cube_flip} , the RSW theory will be used inside~$\calC$. 
Notice that~$\calC(v,w,z,N)$, or rather its infinite version~$\calC(v,w,z,\infty),$ is not biperiodic, 
and therefore, the results of \cite{JiHuangLi, ManoGrimmett} do not apply as such. 
Nevertheless, the RSW theory inside~$\calC(v,w,z,\infty)$ is deduced easily from \cite{JiHuangLi}.

\begin{lemma}[RSW in infinite cubes]
For every~$v,w,z$ as above, ~$\mathcal{C}(v,w,z,\infty)$ satisfies the strong RSW property~${\cal RSW}(c)$ for some~$c > 0$.
\end{lemma}

The constant~$c>0$ above may be shown to be independent of~$v,w,z$, but this is not important for our purposes. 

A piece of notation will be useful: we write~$\mathbb{L}(v,w,z)$ for the isoradial embedding of the square lattice 
formed of one half-plane paved by rhombi of sides~$v,z$ and the other half-plane paved by rhombi of sides~$w,z$; 
the common boundary of the half-planes is the line with direction~$z$ passing through~$0$.

\begin{proof}
Let~$Q$ be a rectangle as in  Definition~\ref{def:RSW} and consider the three pieces~$Q^1,Q^2,Q^3$ of~$Q$ lying inside~$\mathbb{L}(v,w,z)$, 
$\mathbb{L}(w,z,v)$ and~$\mathbb{L}(z,v,w)$ respectively. 
These are finite, potentially empty, subgraphs of square lattices, which have the RSW property by \cite{JiHuangLi}.

Now, deterministically, some~$Q^j$ must contain a translated, rotated and scaled copy of~$[0,4]\times [0,\frac{1}{2}]$, 
and thus, one may cross~$Q$ simply by crossing~$Q^j$. 
This proves the desired by Proposition~\ref{prop:natural_extension} and the RSW property of the associated square lattice.
\end{proof}

For convenience, we will avoid proving Theorem~\ref{thm:cube_flip} for general extensions of the cube and focus on the special case of the infinite cube. 

\begin{lemma}[Reduction to the infinite cube] \label{Lemma:Reduction_to_C}
If Theorem~\ref{thm:cube_flip} holds with~$\mathbb{G}=\calC(v,w,z,\infty),$ then it holds for any natural extension~$\mathbb{G}$.
\end{lemma}

\begin{proof}
Fix~$v,w,z$ and~$N$ as in the statement of Theorem~\ref{thm:cube_flip}. 
Assume Theorem~\ref{thm:cube_flip} holds for~$\bbG = \calC(v,w,z,\infty)$, where we flip only~$ \calC(v,w,z,N)$.
Finally, fix a natural extension~$\bbG$ of~$\calC(v,w,z,\infty)$.

Let~$\omega_1\sim\phi_{\calC(v,w,z,\infty)}$ and~$\omega_2\sim \phi_{(\calC(v,w,z,\infty))^{\mathfrak{F}}}$ (where, in the latter, we only flip~$\calC(v,w,z,N)$)  be coupled so as to satisfy~$(i)$,~$(ii)$, and~$(iii)$ of Theorem~\ref{thm:cube_flip}. 
By Proposition~\ref{prop:natural_extension},~$\omega_1$ may also be coupled to~$\omega\sim \phi_{\mathbb{G}}$ in such a way that
\begin{itemize}
\item $\omega(e)=\omega_1(e)$ for~$e$ in~$\calC(v,w,z,N)$, 
\item for all~$v,w\in \calC(v,w,z,N)$, we have~$v\overset{\omega}{\longleftrightarrow} w$ if and only if~$v\overset{\omega_1}{\longleftrightarrow} w$ and 
\item $\omega$ is conditionally independent of~$\omega_2$ given~$\omega_1$.
 \end{itemize}
Define now a configuration~$\tilde\omega$ on~$\bbG^{\mathfrak{F}}$ by 
$$
\tilde{\omega}(e)=\begin{cases}\omega(e) & e\not\in \calC(v,w,z,N), \\ \omega_2(e) & e\in \calC(v,w,z,N). \end{cases}
$$
Then,~$\tilde{\omega}\sim \phi_{\mathbb{G}^{\mathfrak{F}}}$ by Proposition~\ref{prop:natural_extension} and the fact that the coupling between~$\omega_1$ and~$\omega_2$ satisfies~$(ii)$. Furthermore, the coupling of~$\omega$ and~$\tilde{\omega}$ satisfies~$(i)$,~$(ii)$ and~$(iii)$ because that of~$\omega_1$ and~$\omega_2$ does.
\end{proof}
\subsection{Nail lattices and homotopy distance}

Fix unit vectors~$v,w,z$ not belonging to the same half-plane and~$N\geq1$ large. 
With Lemma~\ref{Lemma:Reduction_to_C} in hand, we now turn to proof of Theorem~\ref{thm:cube_flip} for~$\mathbb{G}=\calC(v,w,z,\infty)$. Write~$\calC = \calC(v,w,z,N)$.

We will encode the distance~$d_{\rm CN}$ via the notion of nails and homotopy classes --- see \cite[Sec. 3.5]{Rotation}. 
We give a brief description below for completeness.   

Several scales will come into play, all depending on~$N$, the size of the cube.
In principle, they are independent, but we will choose them as follows for some small~$\varepsilon > 0$:
\begin{itemize}
\item~$R = N^{1-\varepsilon^3}$ is the distance between the nails; 
\item~$r_{\max} = N^{1-\varepsilon^2}$ is the maximal distance between a nail and the corresponding lattice point; 
\item~$r_{\min} = N^{1-\varepsilon}$ is the scale of mesoscopic clusters;
\item~$r_{\rm IIC} = N^{1-C \varepsilon}$ is the scale at which we will compare extrema of nails with the IIC; $C$ is a large constant chosen as in \cite[eq. (18)]{Rotation}. 
\end{itemize}
Once these quantities are fixed, we define the following notions: We will henceforth assume that $\eps$ is sufficiently small and~$N$ is large, in particular so that~$r_{\rm IIC} \gg 1$.

The lattice~$\mathscr{L}$ is the set of points of~$R\, \bbZ^2$ contained in~$\calC$, at a distance at least~$R$ from~$\calC^c$. 
For a configuration~$\omega$ on~$\bbG$ and a finite cluster~$\sfC$ of~$\omega$, 
define its extremal values in the directions~$-v$ and~$-z$ by 
\begin{align}\label{eq:extrema_def}
	{\rm T}(\sfC) = \max_{y \in \sfC} \langle -v,y\rangle; \quad {\rm B}(\sfC) = \min_{y \in \sfC} \langle -v,y\rangle;\quad   
 {\rm L}(\sfC) = \min_{y \in \sfC} \langle -z,y\rangle;	\quad {\rm R}(\sfC) = \max_{y \in \sfC} \langle -z,y\rangle.
\end{align}
The points realising these minimas and maximas are called the extrema of~$\sfC$. If multiple points satisfy one of the conditions, we choose one according to an arbitrary rule\footnote{The rule is not actually important; in \cite{Rotation} the choice was leftmost, rightmost, bottommost and topmost, respectively, but in retrospect, any choice would do. The actual definition of extremum in \cite{Rotation}  is slightly more complicated, as it involves the extremal cell of the cluster. We will omit this for clarity here.}.

We say that~$\sfC$ is a {\em nail} at some point~$x \in \mathscr{L}$ if 
\begin{align}\label{eq:nail_conditions}
	\big|{\rm T}(\sfC)- \langle -v,x\rangle\big| &\leq r_{\max} ;\quad 
	 &\big|{\rm  B}(\sfC) - \langle -v,x\rangle\big| &\leq r_{\max} ; \nonumber \\  
	\big|{\rm L}(\sfC) - \langle -z,x\rangle\big| &\leq r_{\max} \qquad \text{and}
	&\big|{\rm R}(\sfC) - \langle -z,x\rangle\big| &\leq r_{\max}.
\end{align}

A lattice of nails~$\mathsf N = (\sfN(x))_{x \in \mathscr L}$ is a choice of one nail~$\sfN(x)$ for each lattice point~$x \in \mathscr{L}$. 
Note that there need not be a unique lattice of nails a given configuration - and there might not even be one at all. 

Given a configuration and a lattice of nails~$\sfN = (\sfN(x))_{x \in \mathscr L}$, 
we will encode the large loops and paths of~$\omega \cap \calC$ as follows: 
For each~$x \in \mathscr{L}$, let~$u_x$ be a point in~$\sfN(x)$. For an edge~$xy$ of~$\mathscr{L}$, 
let~$[u_x,u_y]$ be the oriented segment of the straight line between~$u_x$ and~$u_y$.

For a simple oriented curve~$\gamma$ with image in~$\mathbb{R}^2\setminus \{u_x : x \in \mathscr{L}\}$
let~$[u_{x_1}u_{y_1}], \dots, [u_{x_K}u_{y_K}]$ be the sequence of oriented segments crossed by~$\gamma$;
the segments are oriented from left to right as~$\gamma$ crosses them. 
The {\em word} associated to~$\gamma$ is the subsequence of~$(x_1,y_1), \dots, (x_K,y_K)$ 
obtained by inductively removing all consecutive pairs of the form~$(x,y)$,~$(y,x)$. 
When~$\gamma$ is an oriented loop, we consider the word modulo cyclic permutations.
A word and, by extension, the path or loop it encodes is called {\em trivial} if the associated word contains at most one letter.

Recall that the~$\omega$ is encoded by an infinite family of oriented loops.
Write~$[\omega]_{\sfN}$ for the multi-set of non-trivial words of all loops associated to~$\omega$ that intersect $\calC$. 
We call~$[\omega]_{\sfN}$ the homotopy information of~$\omega$ with respect to the lattice of nails~$\sfN$.

The following lemma shows how the homotopy information gives information controls the distance~$d_{\rm CN}$. 

\begin{lemma}\label{lem:nais_to_CN}
	Let	$\bbP$ be a coupling between~$\omega_1\sim \phi_{\mathbb{G}}$ and~$\omega_2\sim \phi_{\mathbb{G}^{\mathfrak{F}}}$ with~$\omega_1|_{\mathcal{C}^c}=\omega_2|_{\mathcal{C}^c}$.
	Then, for any~$\varepsilon > 0$ there exists~$\delta > 0$, independent of $N$, such that 
	\begin{align*}
	&\bbP\big[d_{{\rm global\, CN}}(\omega_1,\omega_2) \leq N^{1-\delta}\big] \\
	&\qquad \leq 	\bbP\big[\text{$\omega_1$,~$\omega_2$ contain lattices of nails~$\sfN_1$,~$\sfN_2$ with~$[\omega_1]_{\sfN_1} = [\omega_2]_{\sfN_2}$}\big] + \tfrac1{N^\delta}.	
	\end{align*}
\end{lemma}

\begin{proof}
	The result is a straightforward extension of  \cite[Prop. 3.9]{Rotation}. 
\end{proof}

\subsection{Coupling for the cube flip}

The cube flip may be viewed as a composition of~$N^2$ finite track-exchanges, where each finite track-exchange is the sequential application of~$N$ star-triangle transformations involving the two same tracks of transverse angles~$v$ and~$w$ and~$N$ tracks of transverse angle~$z$. The track-exchanges will be grouped into~$2N-1$ batches, each batch containing simultaneous, non-overlapping applications of several finite track-exchanges. 

To fix ideas, we will call the rhombi of with sides~$v$ and~$z$ the \emph{old} rhombi (yellow in Figure~\ref{CubePic}), those of sides~$w$ and~$z$ the \emph{new} rhombi (red in Figure~\ref{CubePic}) and the others the \emph{transverse} rhombi (blue in Figure~\ref{CubePic}).
We will perform the cube flip by exchanging old and new rhombi via track-exchanges.
These are described on~$\bbG$, but only affect the region in~$\mathcal{C}$.  We will group them as follows:

Fix~$\bbG_0 = \bbG$. Suppose~$\bbG_t$ was constructed. 
To obtain~$\bbG_{t+1}$, exchange any horizontal track of old rhombi with a horizontal track of new rhombi
if the two tracks are adjacent, and the track of new rhombi is on the top.
Let~$\bfS_t$ denote the finite composition of star-triangle transformations forming these multiple track-exchanges. 

The rhombi in~$\calC$ are gradually permuted, but remain in the same hexagonal region of the plane, which we will always call~$\calC\subseteq \bbG_t$ for clarity. 
In any ~$\bbG_t$,~$\calC$ contains an {\em  active block}  formed of alternating horizontal tracks of old and new rhombi
and up to four {\em frozen blocks}, formed entirely of old, new or transverse rhombi. The latter type of rhombi form two blocks, one to the left and one to the right of the active block. The old and new blocks lie below and above the active block, and exchange positions after~$N$ transformations. 
In~$\bbG_{2N},$ all old rhombi live above all new rhombi, and the process terminates. 
We have~$\bbG_{2N} = \bbG^{\mathfrak{F}}$.

These transformations will also act on configurations to produce a Markov chain of configurations~$\omega_t$ sampled according to~$\phi_{\bbG_t}$. To define the chain, we add a piece of notation. 

A cluster~$\sfC$ in some configuration~$\omega_t$ on some lattice~$\bbG_t$ is called {\em mesoscopic} if 
\begin{align}\label{eq:meso}
	\sfC \subseteq \calC \quad  \text{ and }\quad r_{\min} \leq {\rm T}(\sfC)-{\rm B}(\sfC) &\leq r_{\max}.
\end{align}
Write~${\rm Meso}(\omega_t)$ for the set of mesoscopic clusters of~$\omega_t$ and 
\begin{align}\label{eq:calH_def}
	\calH(\omega) = \big({\rm T}(\sfC), {\rm B}(\sfC), {\rm L}(\sfC), {\rm R}(\sfC)\big)_{\sfC \in {\rm Meso}(\omega)}
\end{align}
for the extremal coordinates of all mesoscopic clusters. 

We now give the explicit construction of the chain~$(\omega_t)_{t = 0, \dots, 2N}$.
Sample~$\omega_0 \sim \phi_{\bbG_0}$. 
For~$t \geq0$ even, given ~$\omega_0,\dots, \omega_t$, we first sample an intermediate configuration~$\omega_{t+1/2} \sim \phi_{\bbG_t}$
as follows:

Consider the boxes~$\La_{r_{\rm IIC}}(z_j)$ around the extrema~$z_j$ of mesoscopic clusters. 
We say that an {\em overlap error} occurs (in~$\omega_t$) if any two such boxes overlap, if any such box intersects more than one cluster of diameter~$r_{\min}/4$, if any such box intersects more than one block, or if any such box is at a distance less than~$r_{\min}$ from~$\calC^c$. In this case, we set~$\omega_{t+1/2} = \omega_t$. 
Otherwise, set~$\omega_{t+1/2} = \omega_t$ for all edges outside of the boxes~$\La_{r_{\rm IIC}}(z_j)$
and resample the configuration in each box according to the measure 
$$\phi_{\bbG_{t}} \Big[ . \,\Big|\, \text{$\omega = \omega_t$ outside~$\bigcup_j\La_{r_{\rm IIC}}(z_j)$ and~$\calH(\omega)=\calH(\omega_t)$} \Big].$$

Now, independently of the definition of~$\omega_{t+1/2}$, set 
\begin{align}
	\omega_{t+1} = \bfS_t(\omega_{t+1/2})\quad \text{ and }\quad \omega_{t+2} = \bfS_{t+1}(\omega_{t}).
\end{align}
The resampling procedure performed between~$\omega_t$ and~$\omega_{t+1/2}$ is described in more detail in \cite[Sec. 3.2]{Rotation}. It is also explained why we resample the environment around the extrema only at even steps.

With these definitions, the main step of the proof is the following statement, corresponding to \cite[Prop 3.7]{Rotation}:
Denote by~${\rm Meso}^+(\omega_0)$ the set of clusters~$\sfC$ of~$\omega_0$ such that 
\begin{align}\label{eq:meso}
	{\rm dist}(\sfC, \calC^c) \geq r_{\max} \quad  \text{ and }\quad 2r_{\min} \leq {\rm T}(\sfC)-{\rm B}(\sfC) &\leq \frac{r_{\max}}2.
\end{align}
For a cluster~$\sfC_{t}$ of~$\omega_t$, we denote by~$\sfC_{t+1}$ its coupled cluster in~$\omega_{t+1}$ with the convention that~$\sfC_{t+1}=\emptyset$ if the cluster ceases to exist. We call a cluster~$\sfC$ of~$\omega_0$ {\em stable} if for all~$t\in [0,2N]$,  we have~$\sfC_t\neq \emptyset$ and furthermore,
\begin{align}
	| {\rm A}(\sfC_0) - {\rm A}(\sfC_t) |  \leq \tfrac{1}{10} r_{\max}
	\qquad \text{ for all~${\rm A} \in \{ {\rm T},{\rm B}, {\rm L}, {\rm R}\}$.}
\end{align}
The above corresponds to \cite[Def. 3.5]{Rotation}.
The constant~$1/10$ above is arbitrary, and should just be taken sufficiently small. 
Note that, as opposed to \cite[Def. 3.5]{Rotation}, we do not mention the drifts~${\rm Drift}_{\rm vert}$ for~$A = {\rm T,B}$ and~${\rm Drift}_{\rm lat}$ in the definition of stability. Indeed, these were proved in  \cite[Lem. 4.4]{Rotation} to be zero.

\begin{prop}\label{prop:stability}
	For~$\eps > 0$ small enough, there exists~$\delta = \delta(\eps) > 0$ 
	so that 
	\begin{align}\label{eq:stability}
		\bbP\big[ \text{all clusters~$\sfC \in {\rm Meso}^+(\omega_0)$ are stable}\big] \geq 1 - \tfrac1{N^\delta}.
	\end{align}
\end{prop}

\begin{proof}
	The proof of \eqref{eq:stability} is the same as that of \cite[Prop 3.7]{Rotation}, with the added input that the drift in both vertical and lateral directions is zero \cite[Lem. 4.4]{Rotation}. 
	
	The only minor difference with \cite[Prop 3.7]{Rotation} is the presence of transfer blocks.  
	Indeed, in  \cite[Prop 3.7]{Rotation}, horizontal tracks are infinite, and no transfer blocks exist; only the old and new blocks do. 
	
	We have incorporated in the notion of overlap errors the situation when 
	an  extremum box intersects the interface between the active and transverse blocks. 
	As such, controlling the expectation of increments for extrema close to such interfaces is not necessary. Also, the probability of an overlap error occurring at a given step may be bounded exactly as 
	in \cite[Prop 3.7]{Rotation}, with the addition of a factor accounting for the case where the overlap error occurs due to the transverse block. 
	
	Another small detail to take into account is in the definition of~$\tau_{\rm start}(\sfC)$ and~$\tau_{\rm end}(\sfC)$ for clusters~$\sfC \in {\rm Meso}^+(\omega_0)$. These are roughly the times when the cluster 
	enters and exits, respectively, the active block. In our setting, a point in the plane may enter or exit the active block by transitioning to the transverse blocks. 
	As such, one needs to take both coordinates of the cluster into account when defining~$\tau_{\rm start}(\sfC)$ and~$\tau_{\rm end}(\sfC)$, rather than just the vertical coordinate, as was the case in  \cite[Prop 3.7]{Rotation}. 
	
	These are only technical difficulties which we will not detail. 
	We will however note that, in this adaptation, the fact that the drifts are null makes all estimates easier. 
\end{proof}

Finally, following \cite{Rotation}, we state a simple consequence of the strong RSW property:

\begin{lemma}\label{lem:nailsmeso}
	For~$\eps > 0$ small enough, there exists~$\delta = \delta(\eps) > 0$ 
	so that 
	\begin{align}\label{eq:nailsmeso}
		\bbP\big[ \text{$\omega_0$ contains a lattice of nails in~${\rm Meso}^+(\omega_0)$}\big] \geq 1 - \tfrac{1}{N^\delta}. 
	\end{align}
\end{lemma}

\begin{proof}
	The proof is identical to \cite[Lem. 3.11]{Rotation}.
\end{proof}

We are now ready to prove Theorem~\ref{thm:cube_flip} for~$\mathbb{G}=\mathcal{C}(v,w,z,\infty)$.

\begin{proof}[Proof of Theorem~\ref{thm:cube_flip} for~$\mathbb{G}=\mathcal{C}(v,w,z,\infty)$]
	Combining Proposition~\ref{prop:stability} and Lemma~\ref{lem:nailsmeso} we find that, under the coupling~$\bbP$ above, for~$\eps >0$ small enough, there exists~$\delta > 0$ such that
	\begin{align*}
	\bbP\big[\text{$\omega_0$ contains a lattice of nails that survives to a lattice of nails of~$\omega_{2N}$ }\big]\geq 1 - \tfrac{2}{N^\delta}. 
	\end{align*}
	
	The construction of the coupling is such that, 
	\begin{itemize}
	\item~$\omega_0 \sim \phi_\bbG$,~$\omega_{2N}  \sim \phi_{\bbG^{\mathfrak{F}}}$, since, by induction, the laws of the different~$\omega_t$ are~$\phi_{\bbG_t}$; 
	\item~$\omega_0 = \omega_{2N}$ outside of~$\calC$; 
	\item all connections between vertices outside~$\calC$ are identical in~$\omega_0$ and~$\omega_{2N}$: indeed, the transformations~$\textbf{S}_t$ preserve such connections,
	and the resampling between~$\omega_t$ and~$\omega_{t+1/2}$ is only performed in situations where no such connection may be broken or created (see the definition of overlap error); 
	\item if~$\sfN_0$ is a lattice of nails in~$\omega_0$ that survives to clusters~$\sfN_{2N}$ in~$\omega_{2N}$, which also form a lattice of nails, then 
	$$ [\omega_0]_{\sfN_0} = [\omega_{2N}]_{\sfN_{2N}}.$$
	\end{itemize} 
	
	Combining the last display above with Lemma~\ref{lem:nais_to_CN} provides the desired conclusion after possibly altering~$\delta$.
\end{proof}

\section{Universality On Bi-Periodic Isoradial Graphs}\label{sec:main_proof}

In this section, we show Theorem~\ref{thm:main}. To be precise, we will prove it in the following form:

\begin{thm}\label{thm:Goalpost}
	Let~$\mathbb{G}$ be a bi-periodic isoradial graph. 
	There exists~$c>0$ such that, for any~$N \geq1$, 
	there exists a coupling~$\bbP$ between~$\omega_1\sim \phi_{\mathbb{G}}$ and~$\omega_2\sim \phi_{\mathbf{Z}^2}$ such that
   $$
    \mathbb{P}\big[{d}_{\rm global\, CN}({\omega_1}|_{\Lambda_N},{\omega_2}|_{\Lambda_N})>N^{1-c}\big]<N^{-c}.
   $$
\end{thm}
In the above, by ${\omega}|_{\Lambda_N},$ we mean the family of loops of $\omega$ contained in $\Lambda_N$. 
The proof of the above occupies the rest of the section. Before diving in, let us mention how it implies  Theorem~\ref{thm:main}.

\begin{proof}[Proof of Theorem~\ref{thm:main} given Theorem~\ref{thm:Goalpost}] 
Fix $\delta > 0$. Let $\bbP$ be the coupling of Theorem~\ref{thm:Goalpost} for $N  = \delta^{-\frac{2}{2-c}}$, where $c > 0$ is the constant given by Theorem~\ref{thm:Goalpost}  --- we will assume for simplicity that $N$ is an integer, otherwise one should round up; we also assume without loss of generality that $c < 1/4$. 

Observe that\footnote{This linearity of the distance is not true for the usual Camia-Newman distance, because rescaling the configuration also changes which loops to consider when determining the distance.}, 
$$
{d}_{\rm global\, CN}(\delta{\omega_1}|_{\Lambda_{\delta N}},\delta{\omega_2}|_{\Lambda_{\delta N}})=\delta \cdot {d}_{\rm global\, CN}({\omega_1}|_{\Lambda_N},{\omega_2}|_{\Lambda_N}).
$$
Then Theorem~\ref{thm:Goalpost} and some standard algebraic manipulations imply
\begin{align}
        \mathbb{P}[{d}_{\rm global\, CN}(\delta {\omega_1}|_{\Lambda_{\delta^{-c/2}}},{\delta \omega_2}|_{\Lambda_{\delta^{-c/2}}})> \delta^{c/4}]
        < \delta^{ c/2},
\end{align}
which in turn yields~\eqref{eq:main_result} for the exponent $c/4$. 
\end{proof}

\subsection{Batched graphs and general strategy}

For an isoradial graph~$\bbG$, we call a family of tracks~$t_1, \dots, t_K$ a~$K$-batch if 
\begin{enumerate}
\item[$i)$] the tracks~$t_j$ do not mutually intersect. Furthermore, they are indexed such that for each~$j$, the tracks~$t_1,\dots,t_{j-1}$ lie on one side of~$t_j$ and~$t_{j+1},\dots, t_K$ on the other;
\item[$ii)$] any rhombus between~$t_1$ and~$t_K$ belongs to exactly one of the tracks~$t_j$. 
\end{enumerate}
We call~$\bbG$ \emph{$K$-batched} if its tracks may be partitioned into~$K$-batches. 
We further call it \emph{properly~$K$-batched} if the tracks of each~$K$-batch have the same transverse angle.

\begin{figure}
\begin{center}
\includegraphics{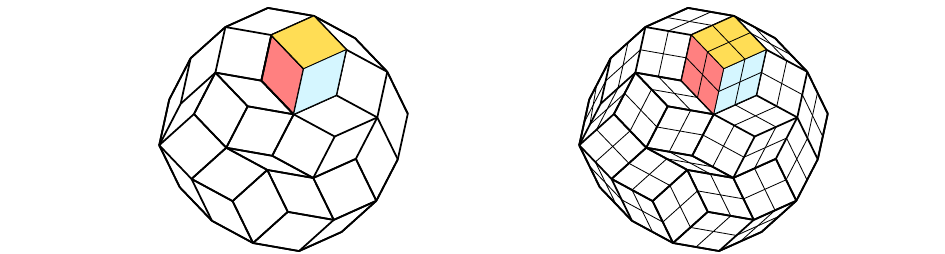}
\end{center}
\caption{An isoradial domain~$\mathcal{D}$ and the corresponding domain in the proper two-batch~$\mathbb{G}^{(2)}$ (scaled by a factor of~$\frac{1}{2}$). Note that stars and triangles of the original rhombic tiling correspond to cubes in~$\mathbb{G}^{(2)}$ to which one may apply the Cube Flip Lemma (Theorem~\ref{thm:cube_flip}).}
\label{fig:two-batch}
\end{figure}

If~$\bbG$ is an isoradial graph, 
its associated proper~$K$-batch, which we write~$\bbG^{(K)}$, is the isoradial graph obtained by replacing each track of~$\bbG$ with~$K$ tracks of the same transverse angle that run next to each other. See Figure~\ref{fig:two-batch}.

Our proof of Theorem~\ref{thm:Goalpost} is split in two steps. First, we will prove that it applies to any properly~$K$-batched bi-periodic isoradial graph, where~$K = N^{1-\delta}$ for some well chosen~$\delta> 0$. This is done by applying Theorem~\ref{thm:cube_flip}. 
We will then prove that if it applies to~$\bbG^{(K)}$, it also does to~$\bbG$. This step requires almost no probabilistic estimates.

\subsection{Proof of Theorem~\ref{thm:Goalpost} for~$K$-batched graphs}

\begin{thm} \label{thm:Kbatched}
Fix a bi-periodic isoradial graph $\bbG$. There exist $C > 0$ and $c \in (0,1)$ such that, 
for any~$K \geq 1$ and $N\geq K$, 
there exists a coupling~$\mathbb{P}$ between~$\omega_1\sim \phi_{{\bbG^{(K)}}}$ and~$\omega_2\sim \phi_{\mathbf{Z}^2}$ such that
    \begin{align}
    \mathbb{P}\left[{d}_{\rm global\, CN}({\omega_1}|_{\Lambda_N},{\omega_2}|_{\Lambda_N})\geq C\left(\tfrac{N}{K}\right)^3 K^{1-c} \right]
    \leq C \left(\tfrac{N}{K}\right)^3 K^{-c}.
    \end{align}
\end{thm}

\begin{remark}\label{rem:Kbatched}
	In the above, $C$ and $c$ are allowed to depend on $\bbG$.
	A quick inspection of the proof shows that they actually only depend on the bounded angles property~\cite{ManoGrimmett} satisfied by $\bbG$. 
	The bi-periodicity of $\bbG$ is only needed here to ensure that $\bbG$ has this bounded angles property. 
	
	We believe that the constants may be taken independently of the bounded angles property, and thus that Theorem~\ref{thm:Kbatched} applies uniformly to all isoradial graphs $\bbG$. This, however, requires controlling cube flips with degenerating angles, which we do not do here. 
\end{remark}

\begin{proof} Fix $N,$ $K$ and a domain $\calD_0$ of $\bbG$ that contains $\Lambda_{N/K}$ and let $\calD$ denote its convexification. 
Then, the perimeter $2M$ of $\calD$ satisfies $M \leq C_0 N/K$ for some constant $C_0$ that depends on $\bbG$. 
Write $e_1,\dots, e_{2M}$ for the edges of the boundary of $\calD$, in counterclockwise order, so that $e_1,\dots, e_M$ have increasing angles in $[0,\pi)$. 
Write $\theta_1,\dots, \theta_M$ for the angles of $e_1,\dots, e_M$. 

Next, add rows of $3M$ rhombi left and right of $\calD$, each with edges with angles $0$ and $\theta_i$ with $i = 1,\dots, M$ except when $\theta_i=0$. This corresponds to prolonging the tracks intersecting $\calD$ horizontally, as in Figures~\ref{fig:completion1} and~\ref{fig:GG}. 
Next, add $3M$ tracks running above the constructed graph, with transverse angle $\pi/2$, and the same below. 
Finally, let $\calG$ be the convexification of the graph thus obtained. 

Note that the edges of the perimeter of $\calG$ have angles:
\begin{itemize}
\item[(I)] $\theta_1,\dots, \theta_M$, each appearing twice; 
\item[(II)] $0$, appearing $12M$ additional times;
\item[(III)] $\pi/2$, appearing an additional $12M$ times.
\end{itemize}
The angles in category (I) different from $0$ and $\frac{\pi}{2}$ are grouped into four sectors, which we denote by $NE$, $NW$, $SW$ and $SE$, respectively. 
The tracks entering $SW$ exit in $NE$ and those entering $NW$ exist in $SW$. 
Also note that the distance between $\Lambda_{N/K}$ and $\calG^c$ is at least $3M$. 

Let $\tilde\calG$ be the paving of the interior of $\calG$ by rhombi constructed as follows: 
\begin{itemize}
\item all tracks entering $SW$ run along the bottom and right side of the boundary up to their exit in $NE$;
\item all tracks entering $NW$ run along the top boundary, then left of the tracks in the point above, up to the time they intersect them and exit in $SE$;
\item all other parts of $\tilde\calG$ are paved by square rhombi - that is, rhombi with angles $0$ and $\pi/2$.
\end{itemize} 
See Figure~\ref{fig:GG} for an example.

\begin{figure}
\begin{center}
\includegraphics[width = 0.4\textwidth, page = 1]{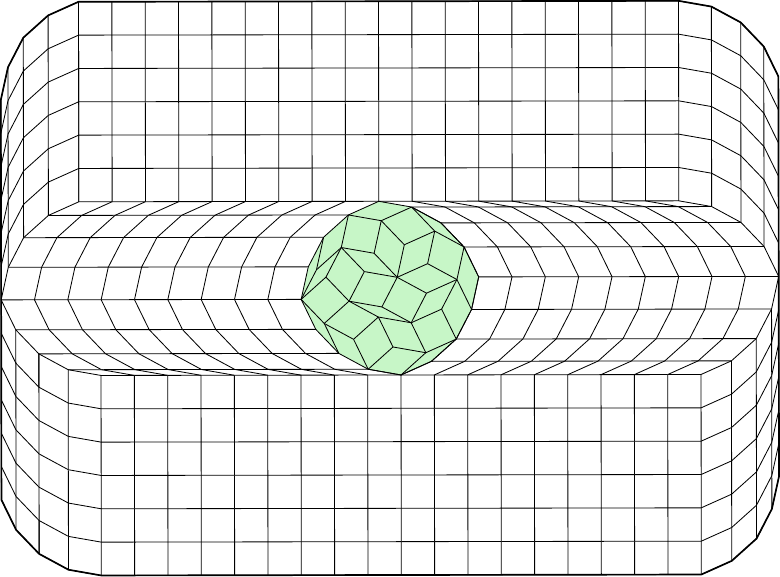}\qquad 
\includegraphics[width = 0.4\textwidth, page = 2]{GG.pdf}
\caption{Left: a convex domain $\calD$ of an isoradial graph (green), with the completion $\calG$ around it.
Right: the convex region occupied by $\calG$, paved with the rhombic tiling $\tilde \calG$. The tracks between SW and NE are blue, those between NW and SE are orange. Note that the green region is entirely covered by a copy of $\bbZ^2$. 
A sequence of at most $C_1 M^3$ star-triangle transformations allows to pass from $\calG$ to $\tilde\calG$.}
\label{fig:GG}
\end{center}
\end{figure}

Observe that, in $\tilde \calG$, $\La_{N/K}$ is entirely contained in the region paved with square rhombi. 
Indeed, the tracks with traverse angles other than $0$ and $\pi/2$ run along the boundary of $\tilde \calG$ and hence, remain within distance $2M$ of it. 

By Lemma~\ref{lem:irred_convex}, there exists a sequence $\bfS_1,\dots, \bfS_L$ of star-triangle transformations with 
\begin{align}
	L \leq C_1M^3 
\end{align} 
that transforms $\calG$ into $\tilde \calG$ where $C_1=13^3$. Here, we have used the counts in (I), (II) and (III) to determine the length of the perimeter.

Let $\bbL$ be a natural extension of $\calG$. 
Then, by applying cube flips $\bfS_1^{(K)},\dots, \bfS_L^{(K)}$ to $\bbL^{(K)}$,
corresponding to the  star-triangle transformations $\bfS_1,\dots, \bfS_L$, 
we transform $\bbL^{(K)}$ into a graph $\tilde \bbL^{(K)}$ which contains $\tilde\calG^{(K)}$. 
Let $\bbP$ be the coupling between $\omega_1 \sim \phi_{\bbL^{(K)}}$ and $\omega_2 \sim \phi_{\tilde \bbL^{(K)}}$ 
obtained by these cube flips. 
By Theorem~\ref{thm:cube_flip}, we find that 
\begin{align}
	\mathbb{P}\big[{d}_{\rm global\, CN}(\omega_1,\omega_2)\geq  C_1M^3 K^{1-c}\big]\leq C_1 M^3\tfrac{1}{K^c}.
\end{align} 
Here, we have used the triangle inequality for the global Camia-Newman distance and a union bound. 

Finally, since $\bbL^{(K)}$ in $\La_{N}$ is identical to $\bbG^{(K)}$,  $\tilde\bbL^{(K)}$ in $\La_{N}$ is identical to $\bfZ^2$,
and since $\bbL^{(K)}$ and $\tilde \bbL^{(K)}$ are natural extensions, Proposition~\ref{prop:natural_extension}
implies that we may extend $\bbP$ to include $\omega_1' \sim \phi_{\bbG^{(K)}}$ and $\omega_2' \sim \phi_{\bfZ^2}$ so that 
\begin{align}
	{d}_{\rm global\, CN}(\omega_1|_{\Lambda_N},\omega_1'|_{\Lambda_N}) = 0 \quad\text{ and } \quad
		{d}_{\rm global\, CN}(\omega_2|_{\Lambda_N},\omega_2'|_{\Lambda_N}) = 0. 
\end{align} 
Combining the last two displays, and using the bound on the value of $M$, we find 
\begin{align}
	\mathbb{P}\big[{d}_{\rm global\, CN}(\omega_1'|_{\Lambda_N},\omega_2'|_{\Lambda_N})\geq  C_2 \big(\tfrac{N}{K}\big)^3 K^{1-c}\big]\leq C_2 \big(\tfrac{N}{K}\big)^3 K^{-c},
\end{align} 
for some $\mathbb{G}$-dependent constant $C_2$, as required. 
\end{proof}

\subsection{Proof of Theorem~\ref{thm:Goalpost} for general graphs}

This section contains the proof of Theorem~\ref{thm:Goalpost}. We start with a proposition. 

\begin{prop}\label{prop:K-batch_univ}
	Fix~$\bbG$ a bi-periodic isoradial graph. 
	For any~$c > 0$, there exists $\delta > 0$ such that, 
	for all~$N \geq 1$, if we set~$K = \lceil N^{1- c}\rceil,$
	there exists a coupling~$\bbP$  between~$\omega_1\sim \phi_{\mathbb{G}}$ and~$\omega_2\sim \phi_{\mathbb{G}^{(K)}}$ such that
   ~$$
    \mathbb{P}\big[{d}_{\rm global\, CN}({\omega_1}|_{\Lambda_N},{\omega_2}|_{\Lambda_N})>N^{1-\delta}\big]<N^{-\delta}.
   ~$$
\end{prop}

\begin{proof}
	We ignore the integer parts for simplicity. 
	We will transform~$\mathbb{G}$ into~$\mathbb{G}^{(K)}$ using star-triangle transformations and track-exchanges. 
	
A collection of tracks~$(t_n)_{n\in \mathbb{Z}}$ is said to be a \textbf{track type} if
\begin{enumerate}
\item[$(i)$] the~$t_n$ are mutually parallel;
\item[$(ii)$] no other track in~$\mathbb{G}$ is parallel to~$t_n$ for any~$n$;
\item[$(iii)$] every other track~$s$ intersects the~$t_n$ in order (i.e. there exists an orientation of~$s$ such that~$t_n$ is intersected before~$t_{n+1}$ for all~$n$). 
 \end{enumerate}
 Since~$\bbG$ is bi-periodic, its tracks may be partitioned into a finite number of track types. 
 Write~$(t^{(j)}_n)_{n\in \mathbb{Z}}$ with~$j = 1,\dots, J$ for these families. \bigskip 
 
\noindent
 We start off with the case when all tracks of the same type have the same transverse angle. 
 Let us start by batching the tracks of the type~$(t^{(1)}_n)_{n\in \mathbb{Z}}$. 
 There exists a finite sequence of transformations~$\bfS_1,\dots, \bfS_T$ with the following properties:
 \begin{itemize}
 	\item each~$\bfS_t$ is the composition of infinitely many, non-overlapping star-triangle transformations;
	\item the tracks~$t^{(1)}_{nK}$ with~$n\in \bbZ$ are not involved in any star-triangle transformation of any~$\bfS_t$;
	\item~$\bfS_T \circ \dots \circ \bfS_1 (\bbG)$ is bi-periodic, with the tracks~$(t^{(1)}_n)_{n\in \bbZ}$ being batched into~$K$-batches. 
 \end{itemize}
The~$\mathbf{S}_t$'s may be constructed by taking limits of the construction in \cite[Lemma 7.1]{ManoGrimmett}. 
That the limit exists is a direct consequence of the fact that the star-triangle transformations are non-overlapping and local. 
 
For a configuration~$\omega_0 \sim \phi_{\bbG}$ set 
\begin{align}
\bbG_1 = \bfS_T \circ \dots \circ \bfS_1 (\bbG) \quad \text{ and }\quad \omega_1 = \bfS_T \circ \dots \circ \bfS_1 (\omega_0).
\end{align}
Note that the connections in~$\omega_0$ and~$\omega_1$ between points on the tracks~$(t^{(1)}_{nK})_{n\in \mathbb{Z}}$ 
are identical. Indeed, the transformations involved in~$\bfS_1,\dots, \bfS_T$ do not affect the points on these track, and thus preserve connections between them. 
Furthermore,~$\omega_1 \sim \phi_{\bbG_1}$. 

Using the same procedure, define sequentially~$\bbG_j$ from~$\bbG_{j-1}$ so 
that each of the families of tracks~$(t^{(1)}_n), \dots , (t^{(j)}_n)$ are batched into~$K$-batches in~$\bbG_j$. 
To get~$\bbG_j$ from~$\bbG_{j-1},$ we apply star-triangle transformations that do not affect any of the tracks~$(t^{(j)}_{nK})_{n\in \bbZ}$. 
Thus, the connections in~$\omega_{j-1}$ and~$\omega_j$ between points on the tracks~$(t^{(j)}_{nK})_{n\in \mathbb{Z}}$ are identical. 

The following lemma is a direct consequence of the RSW property and the preservation of connections between points of $(t^{(j)}_{nK})_{n\in \mathbb{Z}}$. Its proof is postponed to the end of the section.

\begin{lemma} \label{lemma:StripCross}
	For any fixed~$c > 0$, there exists~$\delta > 0$ such that, for all~$N$ large enough and $K$ and $1 \leq j \leq J$ as above, 
   	$$
	\mathbb{P}\big[{d}_{\rm global\, CN}({\omega_{j-1}}|_{\Lambda_N},{\omega_j}|_{\Lambda_N})>N^{1-\delta}\big]<N^{-\delta}.
   	$$
\end{lemma}

We return to the proof of Proposition~\ref{prop:K-batch_univ}.
By repeatedly applying the lemma and using a union bound, we conclude that 
\begin{align}\label{eq:StripCross1}
	\mathbb{P}\big[{d}_{\rm global\, CN}({\omega_{0}}|_{\Lambda_N},{\omega_J}|_{\Lambda_N})> J\cdot N^{1-\delta}\big]< J \cdot N^{-\delta}.
\end{align}
Recall that $\omega_{0} \sim\phi_{\bbG}$ and observe that $\omega_{J} \sim\phi_{\bbG_J}= \phi_{\bbG^{(K)}}$. 
This provides the desired conclusion by adjusting the constants in the case where all track of any given type have the same transverse angle. \smallskip

We now describe how to treat the case where tracks of one type may have  different transverse angles. Due to bi-periodicity, the transverse angles of any track type form a periodic sequence; write $\ell_i$ for the period of the transverse angles of the tracks of type $i$. 
Define $\bbG'$ as the graph obtained by batching $\ell_i K$ tracks of track type $i$ for each $i$. 

The same reasoning as above provides a coupling between $\omega \sim \phi_{\bbG}$ and $\omega' \sim \phi_{\bbG'}$ with 
\begin{align}\label{eq:StripCross2}
	\mathbb{P}\big[{d}_{\rm global\, CN}({\omega}|_{\Lambda_N},{\omega'}|_{\Lambda_N})> C \cdot N^{1-\delta}\big]< C \cdot  N^{-\delta},
\end{align}
where $C$ is a constant that may depend on $\bbG$, but not on $N$ or $c >0$. 
Note, however, that now, $\bbG'$ may be different from $\bbG^{(K)}$, as it may not be properly $K$-batched. 

To get from $\bbG'$ to $\bbG^{(K)}$ we may need to permute tracks of the same type to obtain a properly $K$-batched graph $\bbG''$, 
then move certain batches of $\bbG''$ to obtain $\bbG^{(K)}$. 
More precisely we claim the following: 
There exists a sequence $\bfT_1^{(1)},\dots, \bfT_{L_1}^{(1)}, \dots, \bfT_1^{(J)},\dots, \bfT_{L_J}^{(J)}$, each of which is the composition of infinitely many, non-overlapping track-exchanges
and a sequence $\bfS_1^{(1)}, \dots, \bfS_{T_1}^{(1)}, \dots,  \bfS_1^{(J)},\dots, \bfS_{T_J}^{(J)}$, each of which is the composition of infinitely many non-overlapping star-triangle transformations
such that, 
\begin{align}
( \bfS_{T_J}^{(J)}  \circ \dots \circ \bfS_1^{(1)}\circ \bfT_{L_J}^{(J)}\circ \dots \circ \bfT_1^{(1)})(\bbG') = \bbG^{(K)},
\end{align}
and the transformations $\bfT_j^{(i)}$ and $\bfS_j^{(i)}$ do not affect the tracks $(t^{(i)}_{K\ell_i n})_{n\in \bbZ}$.

In particular, Lemma~\ref{lemma:StripCross} applies to each transformation of the type 
$\bfT_{L_i}^{(i)}\circ \dots\circ\bfT_{1}^{(i)}$ or  $\bfS_{T_i}^{(i)}\circ \dots\circ \bfS_{1}^{(i)}$  and 
extends the coupling $\bbP$ to include a configuration $\omega'' \sim \phi_{\bbG^{(K)}}$ with 
\begin{align}\label{eq:StripCross3}
	\mathbb{P}\big[{d}_{\rm global\, CN}({\omega'}|_{\Lambda_N},{\omega''}|_{\Lambda_N})> C' \cdot N^{1-\delta}\big]< C' \cdot  N^{-\delta},
\end{align}
where $C'$ is a constant that may depend on $\bbG$, but not on $N$ or $c >0$. 
Combining \eqref{eq:StripCross2} and \eqref{eq:StripCross3} provides the desired conclusion after adjustment of the constants. 
\end{proof}

\begin{proof}[Proof of Theorem~\ref{thm:Goalpost}, general case]
We may restrict our attention to $N$ sufficiently large. 
Let $c, C> 0$ be given by Theorem~\ref{thm:Kbatched}. 

Consider $N$ large and let $K = N^{1-c/6}$. 
Applying Theorem~\ref{thm:Kbatched} to $\bbG^{(K)}$ and using Proposition~\ref{prop:K-batch_univ},
we obtain a coupling~$\mathbb{P}$ 
between~$\omega_0\sim \phi_{\bbG}$, 
$\omega_1\sim \phi_{{\bbG^{(K)}}}$ 
and~$\omega_2\sim \phi_{\mathbf{Z}^2}$ such that
    \begin{align}
    &\mathbb{P}\big[{d}_{\rm global\, CN}({\omega_0}|_{\Lambda_N},{\omega_2 }|_{\Lambda_N})> CN^{1 - c/2} + N^{1-\delta} \big]\\
    &\qquad\leq 
	\bbP\big[{d}_{\rm global\, CN}({\omega_0}|_{\Lambda_N},{\omega_1}|_{\Lambda_N})> N^{1-\delta} \big] + 
    \mathbb{P}\big[{d}_{\rm global\, CN}({\omega_1}|_{\Lambda_N},{\omega_2}|_{\Lambda_N})\geq C N^{1- c/2} \big]\\
    &\qquad \leq N^{-\delta} + C N^{-c/3},
    \end{align}
    where $\delta > 0$ is given by Lemma~\ref{lemma:StripCross}. Adjusting the constants provides the desired conclusion. 
\end{proof}

\begin{proof}[Proof of Lemma~\ref{lemma:StripCross}.]
Fix $c > 0$, $N$, $K$ and $j$ as in the statement. We will choose $\delta > 0$ below. We will only work with the tracks of type $j$, which we denote $(t_n)_{n\in \bbZ}$. 

Let~$\mathtt{Strip}_n$ denote the region enclosed by~$t_{nK}$ and~$t_{(n+1)K}$ in~$\Lambda_N$.
Observe that $\La_N$ is then partitioned by $\mathtt{Strip}_{-M}, \dots, \mathtt{Strip}_{M}$ for some $M \leq C N/K = C N^c$, 
where $C$ is a constant potentially depending on $\bbG$. 

Using the fact that connections between the points of $t_{nK}$ are identical in $\omega_{j-1}$ and $\omega_j$, we get a bijection between loops intersecting these tracks in the two configurations.
In order for~${d}_{\rm global\, CN}({\omega_{j-1}}|_{\Lambda_N},{\omega_j}|_{\Lambda_N})>N^{1-\delta}$ to occur, one of two things must happen:
\begin{enumerate}
\item[$a)$] there exists some~$- M \leq n\leq M $ and a loop~$\gamma$ of $\omega_{j-1}$ intersecting~$t_{nK}$
 such that its coupled loop~$\gamma'$ of $\omega_j$ satisfies~$d(\gamma,\gamma')>N^{1-\delta}$ or
\item[$b)$] there exists some~$-M\leq n\leq M$ and a loop~$\gamma$ of $\omega_{j-1}$ or $\omega_j$,  
completely contained in~$\mathtt{Strip}_n$ that has a diameter at least $N^{1-\delta}$. 
\end{enumerate}
Both points imply the existence of a primal path in either $\omega_{j-1}$ or $\omega_j$, 
contained in some $\mathtt{Strip}_n$ with $-M\leq n\leq M$ and of diameter at least $\frac12 N^{1-\delta}$.
Write $\mathtt{Long}$ for this event. 
By the strong RSW property and Proposition~\ref{mixing} applied to $\phi_{\bbG_{j-1}}$ and $\phi_{\bbG_{j}}$ and using a union bound, we find that 
\begin{align}
\bbP[\mathtt{Long} ] \leq C_0 \cdot M \cdot \exp( - c_0 N^{1-\delta}/ K), 
\end{align}
where $C_0$ and $c_0 >0$ are constants potentially depending on $\bbG$. 
Since $K = N^{1-c}$ and $M \leq C N^c$, 
we conclude that, by taking $\delta > 0$ small enough, we have $\bbP[\mathtt{Long} ] \leq N^{-\delta}$. 
This provides the desired conclusion. 
\end{proof}

 \paragraph{Statement on the use of AI,} 
 This work was entirely developed and written without the use of AI. 

 \paragraph{Acknowledgements.} The second author acknowledges financial support from the Swiss National Science Foundation (SNSF).

 \addcontentsline{toc}{section}{References}
\bibliographystyle{plain}
\bibliography{universalitybib2}

@article{Che20,
	author = {Dmitry Chelkak},
	journal = {Annales scientifiques de l'{{\'E}}cole {N}ormale {S}up{\'e}rieure},
	number = {5},
	pages = {1271-1346},
	title = {Ising model and s-embeddings of planar graphs},
	volume = {57},
	year = {2024}}

@article{CheSmi12,
	author = {Chelkak, Dmitry and Smirnov, Stanislav},
	coden = {INVMBH},
	doi = {10.1007/s00222-011-0371-2},
	fjournal = {Inventiones Mathematicae},
	issn = {0020-9910},
	journal = {Invent. Math.},
	mrclass = {82B20 (30C20 30C35 60J67 60K35)},
	mrnumber = {2957303},
	mrreviewer = {Ben Dyhr},
	number = {3},
	pages = {515--580},
	title = {Universality in the 2{D} {I}sing model and conformal invariance of fermionic observables},
	url = {http://dx.doi.org/10.1007/s00222-011-0371-2},
	volume = {189},
	year = {2012}}

@article{KemWer16,
	author = {Kemppainen, Antti and Werner, Wendelin},
	journal = {Probability Theory and Related Fields},
	number = {165},
	pages = {835--866},
	title = {The nested simple conformal loop ensembles in the Riemann sphere},
	year = {2016}}

@article{Ave+26,
	author = {Averous, Emile and Duminil-Copin, Hugo and He, Tiancheng and Lammers, Piet and Manolescu, Ioan},
	note = {preprint arXiv:2610.08507 },
	title = {Conformal invariance of the planar four-state {P}otts model and its random-cluster representation},
	year = {2026}}

@article{PelRom15,
	author = {{R.} Peled and {D.} Romik},
	journal = {S{\'e}minaire Lotharingien de Combinatoire},
	number = {B73b},
	title = {Bijective combinatorial proof of the commutation of transfer matrices in the dense {$O(1)$} loop model},
	volume = {73},
	year = {2015}}

@article{ManoGrimmett,
	author = {Grimmett, Geoffrey R. and Manolescu, Ioan},
	doi = {10.1007/s00440-013-0507-y},
	fjournal = {Probability Theory and Related Fields},
	issn = {0178-8051},
	journal = {Probab. Theory Related Fields},
	mrclass = {60K35 (82B43)},
	mrnumber = {3201923},
	number = {1-2},
	pages = {273--327},
	title = {Bond percolation on isoradial graphs: criticality and universality},
	url = {https://doi.org/10.1007/s00440-013-0507-y},
	volume = {159},
	year = {2014}}

@article{Rotation,
	author = {Duminil-Copin, Hugo and Kozlowski, Karol K. and Krachun, Dmitry and Manolescu, Ioan and Oulamara, Mendes},
	journal = {ArXiV pre-print 2012.11672},
	title = {Rotational invariance in critical planar lattice models},
	year = {2021}}

@article{ADCHLM26,
	author = {Averous, Emile and Duminil-Copin, Hugo and He, Tiancheng and Lammers, Piet and Manolescu, Ioan},
	note = {preprint arXiv:2610.07454},
	title = {Conformal invariance of the six-vertex model's height function},
	year = {2026}}

@article{ScalingRelations,
	author = {Duminil-Copin, Hugo and Manolescu, Ioan},
	doi = {10.1017/fmp.2022.16},
	fjournal = {Forum of Mathematics. Pi},
	issn = {2050-5086},
	journal = {Forum Math. Pi},
	mrclass = {60K35},
	mrnumber = {4515274},
	pages = {Paper No. e23, 83},
	title = {Planar random-cluster model: scaling relations},
	url = {https://doi.org/10.1017/fmp.2022.16},
	volume = {10},
	year = {2022}}

@article{RGHistory,
	author = {Fisher, Michael E.},
	doi = {10.1103/RevModPhys.70.653},
	fjournal = {Reviews of Modern Physics},
	issn = {0034-6861},
	journal = {Rev. Modern Phys.},
	mrclass = {82B28 (82B27)},
	mrnumber = {1627528},
	mrreviewer = {Antonio Gonz\'{a}lez-Arroyo},
	number = {2},
	pages = {653--681},
	title = {Renormalization group theory: its basis and formulation in statistical physics},
	url = {https://doi.org/10.1103/RevModPhys.70.653},
	volume = {70},
	year = {1998}}

@article{SLEBirth,
	author = {Schramm, Oded},
	doi = {10.1007/BF02803524},
	fjournal = {Israel Journal of Mathematics},
	issn = {0021-2172},
	journal = {Israel J. Math.},
	mrclass = {60K35 (30C85 60H15 82B27 82B44)},
	mrnumber = {1776084},
	mrreviewer = {Almut Burchard},
	pages = {221--288},
	title = {Scaling limits of loop-erased random walks and uniform spanning trees},
	url = {https://doi.org/10.1007/BF02803524},
	volume = {118},
	year = {2000}}

@article{RaySpinka,
author = {Ray, Gourab and Spinka, Yinon},
year = { 2020},
title = { A Short Proof of the Discontinuity of Phase Transition in the Planar Random-Cluster Model with $q>4$},
journal  = { Communications in Mathematical Physics},
volume= { 378},
SN  = {1432-0916},
DOi = {10.1007/s00220-020-03827-9}
}

@article{GlazmanLamers,
author = { Glazman, Alexander and Lammers, Piet}, 
year  ={ 2025},
title = {Delocalisation and Continuity in 2D: Loop $\textrm{O}(2)$, Six-Vertex, and Random-Cluster Models},
journal ={ Communications in Mathematical Physics},
volume= {406},
issue ={ 5},
SN  = {1432-0916},
DOI = {10.1007/s00220-025-05259-9}
}

@article{IsingInvariance,
	author = {Chelkak, Dmitry and Duminil-Copin, Hugo and Hongler, Cl\'{e}ment and Kemppainen, Antti and Smirnov, Stanislav},
	doi = {10.1016/j.crma.2013.12.002},
	fjournal = {Comptes Rendus Math\'{e}matique. Acad\'{e}mie des Sciences. Paris},
	issn = {1631-073X},
	journal = {C. R. Math. Acad. Sci. Paris},
	mrclass = {82B20},
	mrnumber = {3151886},
	number = {2},
	pages = {157--161},
	title = {Convergence of {I}sing interfaces to {S}chramm's {SLE} curves},
	url = {https://doi.org/10.1016/j.crma.2013.12.002},
	volume = {352},
	year = {2014}}

@article{DC17,
  title={{Lectures on the Ising and Potts models on the hypercubic lattice.}},
  author={Hugo Duminil-Copin},
  journal={PIMS-CRM Summer School in Probability},
  year={2019},
}

@article{Gri06,
  title= {{The random-cluster model}},
  author={Grimmett, Geoffrey R.},
  journal={volume 333 of { Grundlehren der
	Mathematischen Wissenschaften [Fundamental Principles of Mathematical
	Sciences]}.},
  year={2006},
}

@article{CamiaNewmanFull,
	author = {Camia, Federico and Newman, Charles M.},
	doi = {10.1007/s00220-006-0086-1},
	fjournal = {Communications in Mathematical Physics},
	issn = {0010-3616},
	journal = {Comm. Math. Phys.},
	mrclass = {82B43 (60K35)},
	mrnumber = {2249794},
	mrreviewer = {Christophe Garban},
	number = {1},
	pages = {1--38},
	title = {Two-dimensional critical percolation: the full scaling limit},
	url = {https://doi.org/10.1007/s00220-006-0086-1},
	volume = {268},
	year = {2006}}

@article{TriLat,
	author = {Smirnov, Stanislav},
	doi = {10.1016/S0764-4442(01)01991-7},
	fjournal = {Comptes Rendus de l'Acad\'{e}mie des Sciences. S\'{e}rie I. Math\'{e}matique},
	issn = {0764-4442},
	journal = {C. R. Acad. Sci. Paris S\'{e}r. I Math.},
	mrclass = {60K35 (82B43)},
	mrnumber = {1851632},
	number = {3},
	pages = {239--244},
	title = {Critical percolation in the plane: conformal invariance, {C}ardy's formula, scaling limits},
	url = {https://doi.org/10.1016/S0764-4442(01)01991-7},
	volume = {333},
	year = {2001}}

@article{Dominos,
	author = {Kenyon, Richard},
	doi = {10.1214/aop/1019160260},
	fjournal = {The Annals of Probability},
	issn = {0091-1798},
	journal = {Ann. Probab.},
	mrclass = {52C20 (05B45 05B50 82B41)},
	mrnumber = {1782431},
	mrreviewer = {Mihail N. Kolountzakis},
	number = {2},
	pages = {759--795},
	title = {Conformal invariance of domino tiling},
	url = {https://doi.org/10.1214/aop/1019160260},
	volume = {28},
	year = {2000}}

@article{GFF,
	author = {Schramm, Oded and Sheffield, Scott},
	doi = {10.1007/s11511-009-0034-y},
	fjournal = {Acta Mathematica},
	issn = {0001-5962},
	journal = {Acta Math.},
	mrclass = {60J67 (60D05 60G17 60K35)},
	mrnumber = {2486487},
	mrreviewer = {Julien Dub\'{e}dat},
	number = {1},
	pages = {21--137},
	title = {Contour lines of the two-dimensional discrete {G}aussian free field},
	url = {https://doi.org/10.1007/s11511-009-0034-y},
	volume = {202},
	year = {2009}}

@article{Giuliani1,
	author = {Giuliani, Alessandro and Mastropietro , Vieri and Toninelli, Fabio Lucio},
	doi = {10.1214/15-AIHP710},
	journal = {Annales de l'Institut Henri Poincar{\'e}, Probabilit{\'e}s et Statistiques},
	number = {1},
	pages = {98 -- 168},
	publisher = {Institut Henri Poincar{\'e}},
	title = {Height fluctuations in interacting dimers},
	url = {https://doi.org/10.1214/15-AIHP710},
	volume = {53},
	year = {2017}}

@article{Giuliani2,
	author = {Giuliani, Alessandro and Mastropietro , Vieri and Toninelli, Fabio Lucio},
	doi = {10.1007/s00220-020-03760-x},
	journal = {Commun. Math. Phys.},
	pages = {1883--1959},
	title = {Non-integrable Dimers: Universal Fluctuations of Tilted Height Profiles.},
	volume = {377},
	year = {2020}}

@article{HarmExp,
	author = {Schramm, Oded and Sheffield, Scott},
	fjournal = {The Annals of Probability},
	journal = {Ann. Probab.},
	number = {6},
	page = {2127-2148},
	title = {Harmonic explorer and its convergence to ${S}{L}{E}_4$},
	volume = {33},
	year = {2005}}

@article{BefDum12,
	author = {Beffara, Vincent and Duminil-Copin, Hugo},
	doi = {10.1007/s00440-011-0353-8},
	journal = {Probab. Theory Relat. Fields},
	pages = {511--542},
	title = {The self-dual point of the two-dimensional random-cluster model is critical for $q \geq 1$.},
	volume = {153},
	year = {2012}}

@article{DCcont,
	author = {Duminil-Copin, Hugo and Sidoravicius, Vladas and Tassion, Vincent},
	doi = {10.1007/s00220-016-2759-8},
	fjournal = {Communications in Mathematical Physics},
	issn = {0010-3616},
	journal = {Communications in Mathematical Physics},
	mrclass = {82B26 (82B27 82B43)},
	mrnumber = {3592746},
	mrreviewer = {Farrukh Mukhamedov},
	number = {1},
	pages = {47--107},
	title = {Continuity of the phase transition for planar random-cluster and {P}otts models with {$1 \leq q \leq 4$}},
	url = {https://doi.org/10.1007/s00220-016-2759-8},
	volume = {349},
	year = {2017}}

@article{FKDiscont,
	author = {Duminil-Copin, Hugo and Gagnebin, Maxime and Harel, Matan and Manolescu, Ioan and Tassion, Vincent},
	doi = {10.24033/asens.2485},
	fjournal = {Annales Scientifiques de l'\'{E}cole Normale Sup\'{e}rieure. Quatri\`eme S\'{e}rie},
	issn = {0012-9593,1873-2151},
	journal = {Ann. Sci. \'{E}c. Norm. Sup\'{e}r. (4)},
	mrclass = {82B20 (82B26)},
	mrnumber = {4371394},
	number = {6},
	pages = {1363--1413},
	title = {Discontinuity of the phase transition for the planar random-cluster and {P}otts models with {$q>4$}},
	url = {https://doi.org/10.24033/asens.2485},
	volume = {54},
	year = {2021}}

@article{SuperRSW,
	author = {Duminil-Copin, Hugo and Manolescu, Ioan and Tassion, Vincent},
	doi = {10.1007/s00440-021-01060-6},
	journal = {Probab. Theory Relat. Fields},
	pages = {401--449},
	title = {Planar random-cluster model: fractal properties of the critical phase},
	voume = {181},
	year = {2021}}

@article{JiHuangLi,
	author = {Duminil-Copin, Hugo and Li, Jhih-Huang and Manolescu, Ioan},
	doi = {10.1214/18-EJP223},
	fjournal = {Electronic Journal of Probability},
	journal = {Electron. J. Probab.},
	mrclass = {60K35 (82B20 82B27)},
	mrnumber = {3858924},
	mrreviewer = {Ji\v{r}\'{\i} \v{C}ern\'{y}},
	pages = {Paper No. 96, 70},
	title = {Universality for the random-cluster model on isoradial graphs},
	url = {https://doi.org/10.1214/18-EJP223},
	volume = {23},
	year = {2018}}

@article{electric,
	author = {Kennelly, Arthur E.},
	journal = {Electrical World and Engineer},
	number = 12,
	pages = {413--414},
	title = {The Equivalence of Triangles and Three-Pointed Stars in Conducting Networks},
	volume = 34,
	year = 1899}

@article{Onsager,
	author = {Onsager, Lars},
	fjournal = {Physical Review. Series II},
	issn = {0031-899X},
	journal = {Phys. Rev. (2)},
	mrclass = {82.1X},
	mrnumber = {10315},
	mrreviewer = {J. S. Frame},
	pages = {117--149},
	title = {Crystal statistics. {I}. {A} two-dimensional model with an order-disorder transition},
	volume = {65},
	year = {1944}}

@article{Kenyon,
	author = {Kenyon, Richard},
	doi = {10.1007/s00222-002-0249-4},
	fjournal = {Inventiones Mathematicae},
	issn = {0020-9910},
	journal = {Invent. Math.},
	mrclass = {31C20 (05C50 30G20 82B05)},
	mrnumber = {1933589},
	number = {2},
	pages = {409--439},
	title = {The {L}aplacian and {D}irac operators on critical planar graphs},
	url = {https://doi.org/10.1007/s00222-002-0249-4},
	volume = {150},
	year = {2002}}
\end{document}